\documentclass[reqno,12pt]{amsart}
\usepackage{a4wide}
\usepackage{amscd}
\usepackage{amsmath}
\usepackage{amsfonts}
\usepackage{amssymb}
\usepackage{multirow}
\usepackage{color}
\usepackage{hyperref}

\usepackage[OT2,OT1]{fontenc}
\newcommand\cyr{%
\renewcommand\rmdefault{wncyr}%
\renewcommand\sfdefault{wncyss}%
\renewcommand\encodingdefault{OT2}%
\normalfont
\selectfont}
\DeclareTextFontCommand{\textcyr}{\cyr}

\newcommand{\ca}{{\mathfrak{a}}}

\newcommand{\cg}{{\mathfrak{g}}}

\newcommand{\ck}{{\mathfrak{k}}}

\newcommand{\co}{{\mathfrak{o}}}

\newcommand{\cp}{{\mathfrak{p}}}

\newcommand{\cs}{{\mathfrak{s}}}

\newcommand{\ct}{{\mathfrak{t}}}

\newcommand{\cu}{{\mathfrak{u}}}

\newcommand{\RR}{\mathbb{R}}

\newcommand{\FF}{\mathbb{F}}
\newcommand{\ZZ}{\mathbb{Z}}
\newcommand{\CC}{\mathbb{C}}

\newcommand{\HH}{\mathbb{H}}
\newcommand{\OO}{\mathbb{O}}
\newcommand{\VV}{\mathbb{V}}
\newcommand{\WW}{\mathbb{W}}

\newcommand{\mcD}{\mathcal{D}}

\DeclareMathOperator{\codim}{codim}
\DeclareMathOperator{\Exp}{Exp}

\DeclareMathOperator{\id}{id}
\DeclareMathOperator{\rk}{rk}

\newcommand\bbr{\mathbb{R}}
\newcommand\bbc{\mathbb{C}}

\newtheorem{thm}{Theorem}[section]

\newtheorem{prop}[thm]{Proposition}
\newtheorem{cor}[thm]{Corollary}
\newtheorem{lm}[thm]{Lemma}
\newtheorem{re}[thm]{Remark}

\numberwithin{equation}{section}

\begin {document}

\title{The index of exceptional symmetric spaces}

\author{J\"{u}rgen Berndt}
\address{King's College London, Department of Mathematics, London WC2R 2LS, United Kingdom}
\email{jurgen.berndt@kcl.ac.uk}
\thanks{}

\author{Carlos Olmos}
\address{Facultad de Matem\'atica, Astronom\'ia y F\'isica, Universidad Nacional de C\'ordoba, 
Ciudad Universitaria, 5000 C\'ordoba, Argentina}
\email{olmos@famaf.unc.edu.ar}
\thanks{The second and third author acknowledge financial support from Famaf, UNC and Ciem, CONICET}

\author{Juan Sebasti\'{a}n Rodr\'{i}guez}
\address{Facultad de Matem\'atica, Astronom\'ia y F\'isica, Universidad Nacional de C\'ordoba, 
Ciudad Universitaria, 5000 C\'ordoba, Argentina}
\email{jsrodriguez@famaf.unc.edu.ar }
\thanks{}


\begin{abstract}
The index of a Riemannian symmetric space is the minimal codimension of a proper totally geodesic submanifold (Onishchik \cite{O80}). There is a conjecture by the first two authors (\cite{BO16}) for how to calculate the index. In this paper we give an affirmative answer to this conjecture for the exceptional Riemannian symmetric spaces and for the classical symmetric spaces $Sp_r(\RR)/U_r$. Our methodology is new and based on the idea of using slice representations for studying totally geodesic submanifolds. 
\end{abstract}

\maketitle 

\section {Introduction}

Let $M$ be a connected Riemannian manifold and denote by ${\mathcal S}$ the set of all connected totally geodesic submanifolds $\Sigma$ of $M$ with $\dim(\Sigma) < \dim(M)$. The index $i(M)$ of $M$ is defined by
\[
i(M) = \min\{ \dim(M) - \dim(\Sigma) : \Sigma \in {\mathcal S}\} = \min\{ \codim(\Sigma) : \Sigma \in {\mathcal S}\}.
\]
This terminology was introduced by Onishchik in \cite{O80}. 

Riemannian symmetric spaces are among the most distinguished manifolds in Riemannian geometry. These spaces have been studied by numerous mathematicians and many fascinating properties and applications have been discovered. Despite them being so widely studied, it is a remarkable fact that their totally geodesic submanifolds are not yet known in general. Wolf (\cite{Wo63}) classified them in symmetric spaces of rank $1$ and Klein (\cite{K08},\cite{K09},\cite{K10a},\cite{K10b}) in symmetric spaces of rank $2$.

In previous work (\cite{BO16},\cite{BO17},\cite{BO18}) we developed a systematic approach to the index of Riemannian symmetric spaces. Totally geodesic submanifolds are in one-to-one correspondence to algebraic objects called Lie triple systems. If $M = G/K$ is a Riemannian symmetric space and $\cg = \ck \oplus \cp$ is a corresponding Cartan decomposition, then a Lie triple system is a linear subspace $\VV$ of $\cp$ such that $[[\VV,\VV],\VV] \subseteq \VV$. A distinguished class of Lie triple systems is formed by those Lie triple systems $\VV$ for which the orthogonal complement $\VV^\perp$ of $\VV$ in $\cp$ is also a Lie triple system. Such Lie triple systems are called reflective. Algebraically they correspond to certain involutive automorphisms of $\cg$, geometrically they correspond to totally geodesic submanifolds $\Sigma$ for which the geodesic reflection of $M$ in $\Sigma$ is a well-defined global isometry. These so-called reflective submanifolds have been classified by Leung (\cite{L75},\cite{L79}) in irreducible Riemannian symmetric spaces. Denote by ${\mathcal S}_r$ the set of all connected reflective submanifolds $\Sigma$ of $M$ with $\dim(\Sigma) < \dim(M)$. The reflective index $i_r(M)$ of $M$ is defined by
\[
i_r(M) = \min\{ \dim(M) - \dim(\Sigma) : \Sigma \in {\mathcal S}_r\} = \min\{ \codim(\Sigma) :  \Sigma \in {\mathcal S}_r\}.
\]
It is clear that $i(M) \leq i_r(M)$ and thus $i_r(M)$ is an upper bound for $i(M)$. Moreover, from Leung's work we can calculate $i_r(M)$ explicitly for each irreducible Riemannian symmetric space. This was done explicitly in \cite{BO16}, where we formulated the following conjecture:

\medskip
{\sc Conjecture.} {\it For an irreducible Riemannian symmetric space we have $i(M) = i_r(M)$ if and only if $M \neq G_2^2/SO_4$ and $M \neq G_2/SO_4$. }

\medskip
In our previous work we verified this conjecture for a number of Riemannian symmetric spaces. In \cite{BO17} we verified the conjecture for all exceptional Riemannian symmetric spaces of types II and IV, that is, for the five exceptional compact Lie groups $G_2,F_4,E_6,E_7,E_8$ and their non-compact dual symmetric spaces. The purpose of this paper is to give an affirmative answer to this conjecture for the remaining exceptional Riemannian symmetric spaces, which are of type I (compact) or type III (non-compact). Our main result states:

\begin{thm} \label{mainexc}
Let $M = G/K$ be an irreducible exceptional Riemannian symmetric space. Then $i(M) = i_r(M)$ if and only if $M \neq G_2^2/SO_4$ and $M \neq G_2/SO_4$.
\end{thm}

Duality between symmetric spaces of compact type and of non-compact type preserves totally geodesic submanifolds. We can therefore restrict to symmetric spaces of non-compact type. In Table \ref{exceptionaltable} we list the exceptional Riemannian symmetric spaces of type III together with their index and a totally geodesic submanifold $\Sigma$ with $\codim(\Sigma) = i(M)$.

\begin{table}[ht]
\caption{The index $i(M)$ for irreducible exceptional Riemannian symmetric spaces $M = G/K$ of type III and submanifolds $\Sigma$ of $M$ with $\codim(\Sigma) = i(M)$} 
\label{exceptionaltable} 
{\footnotesize\begin{tabular}{ | p{2.5cm}  p{3cm}  p{1.5cm} p{1cm} p{0.8cm}  p{3cm}   |}
\hline \rule{0pt}{4mm}
\hspace{-1mm}$M$ & $\Sigma$ & $\dim(M)$ & $\rk(M)$ & $i(M)$ & Comments \\[1mm]
\hline \rule{0pt}{4mm}
\hspace{-2mm} 
$E_6^6/Sp_4$ & $F_4^4/Sp_3Sp_1$ & $42$ & $6$ & $14$ &  \\
$E_6^2/SU_6Sp_1$ & $F_4^4/Sp_3Sp_1$ & $40$ & $4$ & $12$ & \\
$E_6^{-14}/Spin_{10}U_1$ & $SO^*_{10}/U_5$ & $32$ & $2$ & $12$ & Onishchik (\cite{O80})\\
$E_6^{-26}/F_4$ & $F_4^{-20}/Spin_9$ & $26$ & $2$ & $10$ & Onishchik (\cite{O80}) \\
$E_7^7/SU_8$ & $\bbr \times E^6_6/Sp_4$ & $70$ & $7$ & $27$ & \\
$E_7^{-5}/SO_{12}Sp_1$ & $E_6^2/SU_6Sp_1$ & $64$ & $4$ & $24$ & \\
$E_7^{-25}/E_6U_1$ & $E_6^{-14}/Spin_{10}U_1$ & $54$ & $3$ & $22$ &  \\
$E_8^8/SO_{16}$ & $\bbr H^2 \times E_7^7/SU_8$ & $128$ & $8$ & $56$ & \\
$E_8^{-24}/E_7Sp_1$ & $E_7^{-5}/SO_{12}Sp_1$ & $112$ & $4$ & $48$ & \\
$F_4^4/Sp_3Sp_1$ & $SO^o_{4,5}/SO_4SO_5$ & $28$ & $4$ & $8$ & Berndt-Olmos (\cite{BO16})\\
$F_4^{-20}/Spin_9$ & $SO^o_{1,8}/SO_8$ $Sp_{1,2}/Sp_1Sp_2$ & $16$ & $1$ & $8$ & Wolf (\cite{Wo63}) \\
$G_2^2/SO_4$ & $SL_3(\RR)/SO_3$ & $8$ & $2$ & $3$ & Onishchik (\cite{O80})\\[1mm]
\hline
\end{tabular}}
\end{table}

The same methodology that we develop in this paper can be used to determine the index of another series of classical symmetric spaces:

\begin{thm} \label{mainsp}
For $M = Sp_r(\RR)/U_r$, $r \geq 3$, we have $i(M) = i_r(M) = 2(r-1)$ and $\Sigma = \RR H^2 \times Sp_{r-1}(\RR)/U_{r-1}$ is a reflective submanifold of $M$ with $\codim(\Sigma) = i(M)$.
\end{thm}

This was proved in \cite{BO16} for $r \in \{3,4,5\}$, but is new for $r > 5$. Taking into account the results from \cite{BO16}, \cite{BO17} and \cite{BO18}, the conjecture remains open for three series of classical symmetric spaces:
\begin{itemize}
\item[(i)] $M = SO^*_{2k+2}/U_{k+1}$ for $k \geq 5$. Conjecture: $i(M) = 2k$.
\item[(ii)] $M = SU^*_{2k+2}/Sp_{k+1}$ for $k \geq 3$. Conjecture: $i(M) = 4k$.
\item[(iii)] $M = Sp_{k,k+l}/Sp_kSp_{k+l}$ for $l \geq 0$ and  $k \geq \max\{3,l +2\}$. Conjecture: $i(M) = 4k$.
\end{itemize}

Our methodology is new and based on the idea of using slice representations for studying totally geodesic submanifolds. In Section \ref{pre} we develop basic sufficient criteria for deciding whether a totally geodesic submanifold is reflective. One of these criteria states that a semisimple totally geodesic submanifold of a symmetric space of rank $\geq 2$ is reflective when the kernel of the full slice representation is non-trivial. It is intuitively clear that fixed vectors of the slice representation must play a crucial rule in the theory. We make this precise in Theorem \ref{main}: a totally geodesic submanifold $\Sigma$ is maximal and the slice representation has a non-zero fixed vector if and only if $\Sigma$ is reflective and the complimentary reflective submanifold is semisimple. This relates to the theory of symmetric $R$-spaces (symmetric real flag manifolds). In Section \ref{fvotsr} we give some applications of Theorem \ref{main} and in Section \ref{extiso} we present its proof. In Section \ref{indxexc} we first derive a very useful lower bound for the codimension of a totally geodesic submanifold. Another important ingredient for our theory is Proposition \ref{hyperbolicfactors}, which states that a reducible totally geodesic submanifold containing a real hyperbolic space or a complex hyperbolic space as a de Rham factor must either be a product of two real hyperbolic spaces or there exists a reflective submanifold of dimension greater than or equal to the dimension of $\Sigma$. This theoretical result will allow us to dismiss many possibilities in our investigations. We eventually apply all these theoretical results to calculate the index of the exceptional symmetric spaces and of $Sp_r(\RR)/U_r$.

\section{Preliminaries and basic results} \label{pre}

Let $M = G/K$ be a connected, locally irreducible, Riemannian symmetric space of compact type, where $G = I^o(M)$ is the identity component of the isometry group $I(M)$ of $M$. We denote by $e$ the identity of $G$ and by $o = eK \in M$ the base point in $M$. Let $\cg = \ck \oplus \cp$ be the corresponding Cartan decomposition of the Lie algebra $\cg$ of $G$. Here, $\ck$ is the Lie algebra of $K$ and $\cp$ is the orthogonal complement of $\ck$ in $\cg$ with respect to the Cartan-Killing form of $\cg$. We denote by $\langle \cdot , \cdot \rangle$ the Riemannian metric on $M$, by $\nabla$ the Riemannian connection on $M$, and by $R$ the Riemannian curvature tensor of $M$. By $r = \rk(M)$ we denote the rank of $M$, which by definition is the maximal dimension of a flat totally geodesic submanifold in $M$. Each $X \in \cg$ determines a Killing vector field $X^*$ on $M$ by $X^*_p = \left.\frac{d}{dt}\right|_{t=0}\Exp(tX)(p)$ for all $p \in M$, where $\Exp : \cg \to G$ is the Lie exponential map. Then $X \in \ck$ if and only if $X^*_o = 0$. It is well known that $[X,Y]^* = -[X^*,Y^*]$ for all $X,Y \in \cg$.

For $p \in M$ we denote by $\sigma_p \in I(M)$ the geodesic symmetry of $M$ at $p$. For $v \in T_oM$ we denote by $\gamma_v : \RR \to M$ the geodesic in $M$ with $\gamma_v(0) = o$ and $\dot{\gamma}_v(0) = v$. Then $\phi_t^v = \sigma_{\gamma_v(t/2)} \circ \sigma_o \in G$ induces a $1$-parameter group $(\phi_t^v)_{t \in \RR}$ of isometries of $M$. Each $\phi _t^v$ is a geometric transvection; it translates the geodesic  $\gamma_v$ by $t$ and $d_{\gamma_v (t_0)} \phi_t^v$ coincides with the parallel transport along $\gamma_v$ from $\gamma_v(t_0)$ to $\gamma_v(t_0+ t)$. If $X^v \in \cg$ with $\phi_t^v = \Exp(tX^v)$, then $X^v \in \cp$, or equivalently $(\nabla X^v)_o = 0$. Conversely, if $X\in \cp$, then $\Exp(tX)$ is a geometric transvection defined by the geodesic $\gamma_v (t) = \Exp (tX)(o)$, $v = X \in \cp \cong T_oM$.
 
Let $\gamma_v$ be a non-trivial closed geodesic in $M$ with $\gamma_v(0) =o$. We may assume, by rescaling $v$, that $\gamma_v$ has (minimal) period $1$. Then $\phi _1^v(o) = o$ and $(\phi^v_1)^2 = \phi ^v_1 \circ \phi ^v_1 = \phi^v _2 = \sigma _{\gamma_v (1)}\circ \sigma _o = \sigma_o^2 = \id _{M}$ is the identity $\id_M$ of $M$. The  isometry $\phi_1^v = \sigma_{\gamma_v(1/2)} \circ \sigma_o$ may be trivial, as is the case when $M$ is a sphere, where $\gamma_v(1/2)$ is the antipodal point of $o$ and $\sigma_{\gamma_v(1/2)} = \sigma _o$. For constructing a non-trivial isometry $\phi_1^v$ on $M$, we have to pass to a suitable globally symmetric quotient of $M$ and then lift back to $M$, provided that $M$ is simply connected. 

Assume that $M$ is simply connected and irreducible. We define an equivalence relation on $M$ by $p \sim q$ if and only if $G_p = G_q$, where $G_p$ and $G_q$ are the isotropy groups of $G$ at $p$ and $q$, respectively. Note that the isotropy groups are connected since $M$ is simply connected and $G = I^o(M)$ is connected. Denote by $\bar{M}$ the quotient space relative to this equivalence relation and by $\pi : M \to \bar{M}$ the canonical projection. Since $M$ is irreducible, the isotropy action of $G_p$ on $T_pM$ is irreducible and therefore has no fixed non-zero vectors. Hence the action of $G_p$ on $M$ has no fixed points apart from $p$ on a sufficiently small open neighborhood of $p$ in $M$. It follows that each fibre of the projection $\pi$ is a discrete subset of $M$. For all isometries $g \in I(M)$ we have $G_{g(p)} = gG_pg^{-1}$ and therefore every isometry $g \in I(M)$ maps equivalence classes to equivalence classes. Thus every $g \in I(M)$ descends to an isometry of $\bar{M}$, where we equip $\bar{M}$ with the induced Riemannian structure from $M$ via $\pi$. We thus can write $\bar{M} = G/G_{\pi(o)}$, where $G$ acts almost effectively on $\bar{M}$. This also shows that every geodesic symmetry of $M$ descends to a geodesic symmetry of $\bar{M}$. Thus $\bar{M}$ is a connected Riemannian symmetric space of compact type. The isotropy group $G_{\pi(o)}$ is not necessarily connected, but its Lie algebra is equal to $\ck$. Note that $\pi : M \to \bar{M}$ induces a bijection from the Lie algebra of Killing vector fields on $M$ onto  the Lie algebra of Killing vector fields on $\bar{M}$.

We denote by $\ck^p \subset \cg$ the Lie algebra of the isotropy group $G_p$ of $G$ at $p \in M$ and by $\ck^{\pi(p)} \subset \cg$  the Lie algebra of the isotropy group $G_{\pi(p)}$ of $G$ at $\pi(p) \in \bar{M}$. Obviously, we have $\ck^p = \ck^{\pi(p)}$. Assume that $\ck^{\pi(p)} = \ck^{\pi(q)}$ for $p,q \in M$. Then we have $\ck^p = \ck^q$. Since the isotropy groups $G_p$ and $G_q$ are connected, this implies $G_p = G_q$ and hence $p \sim q$ and $\pi(p) = \pi(q)$. In other words, different points in $\bar{M}$ have different isotropy algebras. 

\begin{lm} \label{symmetries} Let $\bar{p},\bar{q} \in \bar{M}$ and denote by $\bar\sigma_{\bar{p}},\bar\sigma_{\bar{q}}$ their geodesic symmetries in $\bar{M}$. Then $\bar\sigma_{\bar{p}} = \bar\sigma_{\bar{q}}$ if and only if $\bar{p}=\bar{q}$.
\end {lm} 

\begin{proof}
Let $\bar{s}_{\bar{p}} : G \to G$ be the involutive automorphism of $G$ induced by $\bar\sigma_{\bar{p}}$, that is, $\bar{s}_{\bar{p}}(g) = \bar\sigma_{\bar{p}} g \bar\sigma_{\bar{p}}$ for all $g \in G$. The differential $d_e\bar{s}_{\bar{p}}$ is an involutive automorphism of $\cg$ and $\ck^{\bar{p}}$ is the $+1$-eigenspace of $d_e\bar{s}_{\bar{p}}$. Assume that $\bar\sigma_{\bar{p}} = \bar\sigma_{\bar{q}}$. Then $\bar{s}_{\bar{p}} = \bar{s}_{\bar{q}}$ and therefore $\ck^{\bar{p}} = \ck^{\bar{q}}$, which implies $\bar{p} = \bar{q}$ since two different points in $\bar{M}$ have different isotropy algebras.
\end{proof}

Lemma \ref{symmetries} tells us that $\bar\sigma_{\bar{p}} \neq \bar\sigma_{\bar{q}}$ if ${\bar{p}} \neq {\bar{q}}$. This means that any geodesic symmetry on $\bar{M}$ cannot have another isolated fixed point apart from the obvious one. In other words, there are no poles on $\bar{M}$. The symmetric space $\bar{M}$ is also known in the literature as the adjoint space (\cite{H01}, p.~327) or bottom space (\cite{NT95}, Section 4.2) of $M$. Note that there are simply connected irreducible Riemannian symmetric spaces without poles, that is, $M = \bar{M}$.

Using Lemma \ref{symmetries} we immediately obtain:

\begin {cor}\label{order2} 
Let $\bar{p} \in \bar{M}$ and $\bar{v} \in T_{\bar{p}}\bar{M}$ so that $\gamma_{\bar v}$ is a closed geodesic in $\bar{M}$ with (minimal) period $1$. Then $\bar{q} = \gamma_{\bar{v}}(1/2)$ is an antipodal point of $\bar{p}$ in $\bar{M}$ and $g^{\bar v} = \bar\sigma_{\bar{q}} \circ \bar\sigma_{\bar{p}}$ is a non-trivial involutive isometry of $\bar{M}$ with $g^{\bar v} (\bar{p}) = \bar{p}$. Moreover, $\ell^{\bar v} = d_{\bar{p}} g^{\bar v}$ coincides with the parallel transport in $\bar{M}$ along the geodesic loop $\gamma_{\bar{v} \vert [0,1]}$.
\end {cor}

Since $\bar{q}$ is an antipodal point of $\bar{p}$ in $\bar{M}$, we have $\bar\sigma_{\bar{q}} \circ \bar\sigma_{\bar{p}} = \bar\sigma_{\bar{p}} \circ \bar\sigma_{\bar{q}}$, $\bar\sigma_{\bar{q}}(\bar{p}) = \bar{p} = \bar\sigma_{\bar{p}}(\bar{p})$ and $\bar\sigma_{\bar{p}}(\bar{q}) = \bar{q} = \bar\sigma_{\bar{q}}(\bar{q})$. 

\begin {re} \label{torus} 
\rm In the notation of Corollary \ref{order2}, let $\bar{F}$ be a maximal flat of $\bar M$ with $\bar{p}\in \bar{F}$ and $\bar v\in T_{\bar{p}}\bar{F}$. Such a flat is unique if $\bar{v}$ is a principal vector for the isotropy action of $(G_{\bar{p}})^o$. Since $\bar{F}$ is a torus, it is globally flat. Then, since  $\bar{F}$ is totally geodesic in $\bar{M}$, the restriction $\ell ^{\bar v}\vert_{T_{\bar{p}}\bar{F} } : T_{\bar{p}}\bar{F}  \to T_{\bar{p}}\bar{F}$ is the identity map.  
\end {re}

The quotient space $\bar M$ is only auxiliary and we will lift the isometry $g^{\bar v} \in I(\bar{M})$ to  an isometry $g^v \in I(M)$. Let $p \in M$ with $G_p = K$ (that is, $p \sim o$). Put $\bar{p} = \pi (p) \in \bar{M}$ and choose $v \in T_pM$ so that $\gamma_{\bar v}$ is a closed geodesic in $\bar{M}$ with (minimal) period $1$, where $\bar{v} = d_p\pi(v)$. Then the linear involutive isometry $\ell^v = (d_p\pi)^{-1} \circ \ell ^{\bar v} \circ d_p\pi: T_p M \to T_pM$ preserves the curvature tensor $R^p$ of $M$ at $p$. By the global Cartan Lemma and since $M$ is simply connected, there exists a non-trivial involutive isometry $g^v$ of $M$ with $g^v(p) = p$ and $d_p g^v = \ell ^v$. The isometry $g^v \in I(M)$ is not necessarily in $G = I^o(M)$. Note that the linear isometry $\ell^v$ also induces an isometry of the dual symmetric space of $M$, since the curvature tensor of the dual is just $-R$.

Note that any flat $\bar{F}$ of $\bar{M}$ can be written as $\bar{F} = \pi (F)$ with some flat $F$ of $M$. 

\begin {prop}\label{reflective} 
Let $M = G/K$ be a simply connected, irreducible, Riemannian symmetric space with $G= I^o(M)$,  $o= eK$ and $\cg = \ck \oplus \cp$ the Cartan decompositon of $\cg$ at $o$.  Let $0\neq w\in T_oM \cong \cp$ be arbitrary.  Then there exists a non-trivial isometry $g$ of $M$ with the following properties:
\begin{itemize}
\item[(i)] $g(o) = o$;
\item[(ii)] $g^2 = {\rm id}_M$;
\item[(iii)] The orbit $K \cdot w \subset \cp$ is invariant under $d_og$, that is, $d _og (K\cdot w) = K\cdot w$;
\item[(iv)] Each vector in $\nu _w (K\cdot w) \cong Z_\cp(w)$ is fixed by $d _og$, where $\nu _w (K\cdot w)$ is the normal space of $K \cdot w$  at $w$ and $Z_\cp(w)$ is the centralizer of $w$ in $\cp \cong T_oM$; 
\item[(v)] Assume that  $\rk(M) \geq 2$, or equivalently, that $K\cdot w$ is not a sphere.  Then the linear subspace $\nu _w (K\cdot w)$ of $T_oM$ coincides with the set of fixed vectors of $d_og$ if and only if $K\cdot w$ is an extrinsically symmetric orbit. 
\end{itemize}
\end{prop}

\begin {proof} 
It suffices to prove the proposition for the case that $M$ is of compact type. 

We have $Z_\cp(w) = \bigcup_{\ca \in {\mathcal A}^w} \ca$, where ${\mathcal A}^w$ is the set of maximal abelian subspaces $\ca$ of $\cp$ with $w \in \ca$. The intersection $\ca_0 = \bigcap_{\ca \in {\mathcal A}^w} \ca$ is the abelian part of $Z_\cp(w)$. Moreover, as it is standard to prove, there exists a compact flat $F$ of $M$ with $o \in F$ and $T_oF = \ca_0$.  Now write $K \cdot w = K/K_w$, where $K_w$ is the isotropy group of $K$ at $w$, and let $K_w^o$ be the identity component of $K_w$. It is clear that $\ca_0$ is invariant under the action of $K_w^o$ on $\cp$, and therefore $F$ is invariant under the action of $K_w^o$ on $M$.
Since the connected group $K_w^o$ fixes $o$, we obtain that $K_w^o$ acts trivially on $F$. It follows that $K_w^o$ acts trivially, via the isotropy representation, on the subspace $\ca _0$ of $Z_\cp(w) = \nu _w (K\cdot w) $. In other words, the connected slice representation of $K_w$ on the normal space $\nu _w (K\cdot w)$ fixes $\ca_0$ pointwise. 

Now consider the bottom space $\bar{M}$ and the canonical projection $\pi : M \to \bar{M}$. The image $\bar{F} = \pi(F)$ is a compact flat of $\bar{M}$. Choose $0 \neq v\in \ca _0 =  T_oF$ such that the geodesic $\gamma _{\bar{v}} = \pi \circ \gamma _v$, $\bar{v} = d_o \pi (v)$, is closed in $\bar{F}$. Since the initial directions of closed geodesics in $\bar{F}$ starting at $\bar{o} = \pi(o)$ form a dense set in $T_{\bar{o}}\bar{F} \cong T_oF = \ca_0$, we can choose $v$ arbitrarily close to $w$. The orbits $K\cdot w$ and $K\cdot v$ are parallel orbits if $v$ is sufficiently close to $w$  (see \cite{BCO16}, Corollary 2.3.7). For such $v$ this then implies $\nu _v (K\cdot v) = \nu _w (K\cdot w)$, or equivalently, $Z_\cp(v) = Z_\cp(w)$. We normalize $v$ so that $\gamma_{\bar v}$ is a closed geodesic in $\bar{M}$ with (minimal) period $1$ and construct the isometries $g^{\bar{v}} \in I(\bar{M})$ and $g^v \in I(M)$ as above. Using Remark \ref{torus}, and using the fact that $Z_\cp(\bar v) = Z_\cp(v)$ is the union of all abelian subspaces of $\cp \cong T_{\bar{o}}\bar M$ containing $\bar v$, we obtain that $d_{\bar{o}} g^{\bar v}$ fixes $Z_\cp(\bar v) = Z_\cp(v)$ pointwise and, in particular, $d_{\bar{o}} g^{\bar v}(\bar v) = \bar v$. Then $g = g^v \in I(M)$ satisfies the properties (i)--(iv) (for (iii) use that $gKg^{-1} = K$, since $K$ is the identity component of the full isotropy group of $M$).  
 
The ``only if'' part of (v) is just the definition of an extrinsically symmetric space. It remains to prove that,  if $S= K\cdot w$ is an extrinsically symmetric space, then $d_og$ coincides with the extrinsic symmetry (and so it has no fixed vector tangent to the orbit $K\cdot w$). So let us assume that $S$ is extrinsically symmetric. It is well known that in this case the orbit $K\cdot w$ must be most singular, that is, nearby orbits in the sphere have greater dimension than $K \cdot w$, or equivalently, $\dim(\ca_0) = 1$ and so $v$ is a scalar multiple of $w$, hence we may assume that $v=w$ (note that this also follows from Theorem 4.2 in \cite{BO16}). Let 
\[
\tilde  K = \{ h\in SO (T_oM): h (K\cdot v)=  K\cdot v\}^o.
\]
Note that $\tilde K$ contains the transvections $\rho _p \circ \rho  _q$, where $\rho_p,\rho_q$ are the extrinsic symmetries at $p,q \in K\cdot v$. Then $\tilde K \supseteq K$ and $\tilde K$ is not transitive on the sphere of $T_pM$, since $\tilde K\cdot v = K\cdot v$. From Simons Holonomy Theorem  (\cite{S62}) we then get $\tilde K = K$ (see Remark 8.3.5 in \cite {BCO16}) and it follows that $(K,K_v)$ is a symmetric pair.  
 
Since $K_v \cdot v = \{v\}$, the construction of $g= g^v$ shows that $k g k^{-1}= g$ for all $k\in K_v$. Let $\VV = \{u \in T_v(K\cdot v) : d_og(u) = u\}$. Then $\VV$ extends to a $K$-invariant parallel distribution $\mcD$ on the orbit $S= K \cdot v$.  The orthogonal complement $\VV^\perp$ of $\VV$ in $T_vS$  is exactly the $(-1)$-eigenspace 
of the linear isometry  $d_og$. Let $\alpha$ be the second fundamental form of $S$. 
Then we have
\[
\alpha (u,z) = d_og (\alpha (u,z)) = \alpha (d_og(u),d_og(z)) = \alpha (u, -z) = - \alpha (u,z)
\]
for all $u\in \mathbb V$ and $z\in \mathbb V^\perp$, and hence $\alpha (\mathbb V, \mathbb V^\perp) = 0$.
This implies $\alpha (\mcD , \mcD^\perp)=0$. If $\mcD$ is non-trivial, Moore's Lemma tells us that $S$ splits extrinsically (see \cite {BCO16}, Lemma 1.7.1). In this case $K$ does not act irreducibly, which is a contradiction. We conclude that $d_og$ coincides with the extrinsic symmetry $\rho_v$ of $S$.   
\end{proof}

\begin{re} \rm  The above construction of the symmetry for extrinsically symmetric spaces is due to Nagano in the context of fibrations of symmetric spaces (see \cite{NT95} and its bibliography; see also Section 7 in \cite{E17}).  However,  in Nagano's construction the closed geodesic must be minimizing until it reaches its antipodal point. In the proof of Proposition \ref{reflective}, the closed geodesic does not generally have this property. Let $M= G/K$ be a simply connected irreducible symmetric space. If $S= K\cdot v$ is a most singular isotropy orbit, which is not extrinsically symmetric, then $\exp_o (tv)$ is a closed geodesic that is not minimizing until it reaches its antipodal point. Namely, if $\alpha _1 ,\ldots , \alpha _r$ is a basis of simple roots, then $K\cdot v = K\cdot H^i$ for some $i$, 
where $H^1 ,\ldots , H^r$ is the dual basis. From Kobayashi and Nagano \cite{KN64}, $S$ is extrinsically symmetric if and only if $\delta _i = 1$,  where  $\delta = \delta _1 \alpha_1 + \ldots + \delta_r \alpha_r$ is the highest root. 
From  \cite{H01}, Chapter VII-3, the closed geodesic $\gamma_{H^i}$ is not minimizing beyond $\frac{1} {2\delta_i}$ of its length. 
\end {re}

Let $M= G/K$ be a simply connected irreducible symmetric space and $\Sigma$ be a non-semisimple maximal totally geodesic submanifold of $M$ with $o \in M$. Then $T_o\Sigma = Z_\cp(v)$, where $K \cdot v$ is a most singular orbit. Note that,  if $K \cdot v$ is a most singular orbit, then $Z_\cp(v)$ is a non-semisimple Lie triple system of  $T_oM \simeq \cp$ whose abelian part is $\RR v$. It was proved in \cite{BO16} (Theorem 4.2) that $Z_\cp(v)$, which coincides with the normal space  $\nu _v(K\cdot v)$,  is the tangent space to a maximal non-semisimple totally geodesic submanifold of $M$ if and only if $K\cdot v$ is an extrinsically symmetric orbit. If $K\cdot v$ is extrinsically symmetric, then both $\nu _v(K\cdot v)$ and $T_v(K\cdot v)$ are Lie triple systems in $T_oM$. In fact, $T_v(K\cdot v)$ coincides with the fixed vectors of $d_o\sigma_o  \circ d_og^v$, where $\sigma_o$ is the geodesic symmetry of $M$ at $o$ and $d_og^v$ is the extrinsic symmetry of $K\cdot v$ at $v$ (with $g^v$ constructed as above).

Let $\Sigma$ be a connected, complete, totally geodesic submanifold of $M$ with $o \in \Sigma$. Then $\Sigma$ is called {\it reflective} if the normal space $\nu_o\Sigma$ is a Lie triple system as well. Note that  $\Sigma $ is reflective if and only if $\Sigma$ is the connected component containing $o$ of the 
fixed point set of an involutive isometry $\tau \in I(M)$ with $\tau(o) = o$ (see, for instance, \cite{BO16}). 

In the proof of \cite{BO16}, Theorem 4.2, it is shown that $Z_\cp(v)$ is properly contained in a (proper) Lie triple system of $T_oM$ if $K\cdot v$ is not extrinsically symmetric. The following corollary of Proposition \ref {reflective} gives an alternative proof of this result with   additional information: $Z_\cp(v)$ is properly contained in a (proper)  reflective Lie triple system of $T_oM$. 

\begin {cor}
Let $M = G/K$ be a simply connected irreducible Riemannian symmetric space and $0 \neq w \in T_oM$. 
Assume that $K\cdot w$ is not  extrinsically symmetric. Then $Z_\cp(w)$ is properly contained in a reflective (proper) Lie triple system of $T_oM$. 
\end{cor}

\begin {proof} 
Let $g$ be as in Proposition \ref {reflective}. Since $g$ is of order $2$, the set of fixed vectors of $d_og$ is a reflective Lie triple system containing $Z_\cp(w)$. 
\end {proof}

We finish this section by proving that if a connected, complete, totally geodesic submanifold $\Sigma$ of a symmetric space $M$ contains a reflective submanifold of $M$, then $\Sigma$ must be reflective as well. For this aim we first generalize Proposition 3.4 of \cite{BO16} to include the case that the kernel of the slice representation is finite. 

Let  $\Sigma = G^\Sigma/K^\Sigma$ be a semisimple totally geodesic submanifold of a (simply connected) symmetric space $M=G/K$ with $o \in \Sigma$, where 
\[
G^\Sigma = \{g\in I (M) : g(\Sigma) = \Sigma\}, 
\]
$K^\Sigma = (G^\Sigma)_o $ is the isotropy group of $G^\Sigma$ at $o$ and $I(M)$ is the full isometry group of $M$ (whose identity component coincides with $G$). The group $G^\Sigma$ is, in general, neither connected nor effective on $\Sigma$. Note that $G^\Sigma$ contains the \textit{glide transformations} of $\Sigma$, that is, the closed subgroup of $G^\Sigma$ with Lie algebra $[T_o\Sigma,T_o\Sigma] \oplus T_o\Sigma$. 

The \textit{full slice representation} of $\Sigma$ at $o$ is the representation $\rho : K^\Sigma \to O(\nu _o\Sigma)$ given by $\rho(k)= d_ok_{\vert \nu_o \Sigma}$, where $\nu _o\Sigma = (T_o\Sigma)^\perp$ is the normal space of $\Sigma$ at $o$.

\begin {re}\label {dualsigma} \rm If $M = G/K$ and $\hat{M} = \hat{G}/K$ are dual symmetric spaces, then $\cg = \ck + \cp \subset \cg^\CC$ and $\hat{\cg} = \ck + i\cp \subset \cg^\CC$ are the corresponding Cartan decompositions at $o = eK \in M$ and $\hat{o} = eK \in \hat{M}$. The isotropy representation of $G/K$ at $o$ on $T_oM \cong \cp$ is canonically equivalent to the isotropy representation of $\hat{G}/K$ at $\hat{o}$ on $T_{\hat{o}}\hat{M} \cong i\cp$. If $V \subset \cp$ is a Lie triple system in $\cp$, then $iV \subset i\cp$ is a Lie triple system in $i\cp$, and vice versa. Thus we have a natural bijection between the Lie triple systems on $\cp$ and $i\cp$, and therefore between connected, complete, totally geodesic submanifolds in $M$ and $\hat{M}$. 

Let $\Sigma$ be a (connected, complete) totally geodesic submanifold of $M$ with $o\in \Sigma$ and let $\hat{\Sigma}$ be the (connected, complete) totally geodesic submanifold of $\hat{M}$ with $T_{\hat{o}}\hat{\Sigma} = iT_o\Sigma$. Then $K^\Sigma = K^{\hat {\Sigma}}$ (via the respective isotropy representations). In fact, if $k\in K^\Sigma$, then $-i(d_ok) i$ is a linear isometry of $i\cp = T_{\hat o} \hat{M}$ leaving  $iT_o\Sigma$ invariant and preserving the curvature tensor $\hat{R}^{\hat o} $ of $\hat{M}$ at $\hat{o}$, since $\hat{R}^{\hat o}(ix,iy)iz= -iR^o(x,y)z$ for all $x,y,z \in T_oM$. So, by the global Cartan Lemma, $-i(d_ok)i$ is the differential of an isometry of $\hat M$ fixing $\hat{o}$. In particular, $\Sigma$ is reflective if and only if $\hat{\Sigma}$ is reflective. 
\end {re}

\begin {prop}\label {fullslice} 
Let $M =G/K$ be an irreducible Riemannian symmetric space with $\rk(M) \geq 2$, where $(G,K)$ is an effective Riemannian symmetric pair.  Let $\Sigma$ be a semisimple totally geodesic submanifold of $M$ with $o \in \Sigma$. Assume that the kernel of the full slice representation of $\Sigma$ at $o$ is non-trivial. Then $\Sigma$ is reflective.
\end {prop}

\begin {proof} 
The proof is similar to the one of Proposition 3.4 in \cite{BO16}. By Remark \ref{dualsigma} we may assume that $M$ is of non-compact type and that $\Sigma $ is complete (and therefore diffeomorphic to $T_o\Sigma$). 

Let $K^\Sigma = \{k \in I(M)_o: k(\Sigma) =  \Sigma\}$ and $\rho : K^\Sigma \to O(\nu _o\Sigma)$ be the full slice representation. The kernel $H$ of $\rho$ is a normal subgroup of $K^\Sigma$ (note that $H$ may be finite).  The subspace 
\[
\VV = \{ v \in T_o\Sigma : d_ok(v) = v \mbox{ for all } k \in H\} 
\]
is a Lie triple system in $T_o\Sigma$. Moreover, since $H$ is a normal subgroup of $K^\Sigma$, the subspace $\VV \subseteq T_o\Sigma$ is invariant under $K^\Sigma$ and, in particular,  under the isotropy group $K^\prime = (G')_o$ of the group $G^\prime$ of glide transformations of $\Sigma$. Thus $\mathbb V$ extends to a $G'$-invariant, and hence parallel, distribution on $\Sigma$. It follows that $\Sigma$ is the Riemannian product $\Sigma = \Sigma _1\times \Sigma _2$ of two totally geodesic submanifolds $\Sigma_1$ and $\Sigma_2$ with $T_o\Sigma _1= \VV$ and $T_o\Sigma _2= \VV^\perp \cap T_o\Sigma$. 

The subspace $\VV \oplus \nu_o\Sigma$ is a Lie triple system in $T_oM$, since it consists of the fixed vectors of the action of $H$ on $T_oM$. This implies that $\Sigma_2$ is reflective. We write $\Sigma _2 = G''/K''$, where $G'' \subseteq G$ is the group of  glide transformations of $\Sigma _2$ ($G''$ acts almost effectively on $\Sigma _2$) and $K'' = (G'')_o$. Since $\Sigma = \Sigma _1\times \Sigma _2$, the group $K''$ acts trivially on $T_o\Sigma _1 =\VV$. Consequently $\VV$ is a subset of the set $\WW$ of fixed vectors of the action of $K''$ on $T_o \Sigma _3$, where $\Sigma _3$ is the totally geodesic submanifold of $M$ associated with the Lie triple system $\VV \oplus \nu_o\Sigma$. By  Lemma 3.2 (ii) in \cite{BO16}, the subspace $\WW$ is the tangent space of a totally geodesic flat submanifold $\Sigma _0 \subseteq \Sigma _3$. Since $T_p\Sigma _1 = \VV \subseteq \WW$, we conclude that $\Sigma _1$ is flat, which is a contradiction since $\Sigma$ is semisimple unless $ \VV =\{0\}$. It follows that $\Sigma = \Sigma_2$ is reflective. 
\end {proof}

\begin {cor}  \label  {proper-reflective}
Let $M =G/K$ be an irreducible  simply connected Riemannian symmetric space with $\rk(M)\geq2$. Let $\Sigma_1,\Sigma _2$ be connected, complete, totally geodesic submanifolds of $M$ with $\Sigma_1 \subseteq \Sigma_2$. If $\Sigma _1$ is reflective, then $\Sigma_2$ is reflective.
\end {cor}

\begin {proof}
We can assume that $M$ is of non-compact type and $o \in \Sigma_1$. 

If $\Sigma _1$ is non-semisimple, it follows from Proposition 4.2 in \cite{BO18} that $T_o\Sigma _1$ is the normal space of  a symmetric isotropy orbit of $K$. Then, by Theorem 1.2 of \cite{BO16}, $\Sigma _1$ is maximal and hence $\Sigma_2 = \Sigma_1$ is reflective.

If $\Sigma _1$ is semisimple, we consider two cases:

Case 1: $\Sigma _2$ is semisimple. Let $\tau$ be the geodesic reflection of  $M$ in the reflective submanifolds $\Sigma _1$ and $\sigma $ be the geodesic symmetry of $M$ in $o$. Then $h= \sigma  \circ \tau$ is involutive and the eigenspaces of $d_oh$ are $T_o\Sigma _1$ (associated with the eigenvalue $-1$) and $\nu_o\Sigma_1$ (associated with the eigenvalue $1$). We decompose $T_o\Sigma_2$ into  $T_o\Sigma _2 = T_o\Sigma _1 \oplus \VV$ with $\VV\subseteq \nu_o \Sigma _1$. Then $d_oh(T_o\Sigma _2) = T_o\Sigma_2$ and thus $h \in K^{\Sigma _2}$. Moreover, the non-trivial isometry $h$  belongs to the kernel of the full slice representation of $\Sigma _2$ at $o$. It follows from Proposition \ref{fullslice} that $\Sigma _2$ is reflective.

Case 2: $\Sigma _2$ is not semisimple. Write $\Sigma _2 = \Sigma _0 \times \Sigma _s$ (Riemannian product), where $\Sigma_0$ is the Euclidean factor of $\Sigma_2$ and $\Sigma_s$ is the semisimple factor of $\Sigma_2$. Obviously, we have $\Sigma _1 \subseteq \Sigma _s$. From Case 1 we conclude that $\Sigma_s$ is reflective.
However, the semisimple  part of a non-semisimple totally geodesic submanifold is never reflective, due to Corollary 3.3 in \cite{BO16}. So Case 2 cannot occur. 
\end {proof}

\section{Fixed vectors of the slice representation} \label{fvotsr}

Let $M=G/K$ be an irreducible, simply connected, Riemannian symmetric space with $\rk(M) \geq 2$, where $G = I^o(M)$, $K=G_o$ and $o \in M$. The corresponding Cartan decomposition at $o$ is $\cg = \ck \oplus \cp$ and the tangent space $T_oM$ is identified with $\cp$ in the usual way. 

Let $\Sigma$ be a complete totally geodesic submanifold of $M$ with $o \in \Sigma$ and put $T_o\Sigma \cong \cp' \subseteq \cp$. Let $G' \subseteq G$ be the subgroup of $G$ consisting of the glide transformations of $\Sigma$ and let $K' = (G')_o$. The Lie algebras of $K'$ and $G'$ are $\ck' =[\cp',\cp'] \subseteq \ck$ and $\cg'  = \ck' \oplus \cp' \subseteq \ck \oplus \cp = \cg$ respectively. Then $\Sigma = G'/K'$ and $(G',K')$ is an almost effective symmetric pair. Note that normal vectors of $\Sigma$ that are fixed by the slice representation correspond to $G$-invariant normal vector fields that are parallel along $\Sigma$.
  
We now state one of our main results; the proof will be given in the next section.
  
\begin {thm} \label{main}  
Let $M=G/K$ be an irreducible, simply connected, Riemannian symmetric space with $\rk(M) \geq 2$, where $G = I^o(M)$, $K=G_o$ and $o \in M$. Let $\Sigma = G'/K'$ be a (proper) totally geodesic submanifold of $M$ with $o \in \Sigma$ and $\dim(\Sigma) \geq \frac{1}{2}\dim(M)$, where $G'$ is the subgroup of $G$ consisting of the glide transformations of $\Sigma$ and $K' = (G')_o$. Let $\rho$ be the slice representation of $(K')^o$ on $SO(\nu_o\Sigma)$. Then the following statements are equivalent: 
\begin{itemize}
\item[(i)] $\Sigma$ is maximal and there exists a non-zero vector in $\nu_o\Sigma$ that is fixed by the slice representation $\rho$. 
\item[(ii)] $\Sigma$ is reflective and the complementary reflective submanifold is non-semisimple.
\item[(iii)] $T_o\Sigma$ coincides, as a linear subspace, with the tangent space $T_v(K \cdot v)$ of a symmetric isotropy orbit.
\end{itemize}
\end {thm}  

Note that the assumption in (i) for $\Sigma$ to be maximal is necessary. Consider for example the symmetric space $Sp_r /U_r$ and its maximal totally geodesic submanifold $Sp_1/U_1 \times Sp_{r-1}/U_{r-1} = \RR H^2 \times Sp_{r-1}/U_{r-1}$. If we now consider $\Sigma = Sp_{r-1}/U_{r-1}$, then the slice representation of $U_{r-1}$ fixes any vector in $T_o \RR H^2$.
  
It follows from the classification of symmetric $R$-spaces (see \cite{KN64}) that a non-semisimple extrinsically symmetric isotropy orbit (briefly, symmetric isotropy orbit) $K \cdot v \subset T_oM$ has always half the dimension of $T_oM$ and that the Euclidean local factor of $K \cdot v$ is $1$-dimensional. Moreover,  $T_v(K\cdot v)$ and $\nu_v(K\cdot v)$ are $K$-equivalent Lie triple systems. 
In the next proposition we give a conceptual proof of these observations. 

\begin {prop} \label{half} 
Let $M = G/K$ be an irreducible, simply connected, Riemannian symmetric space with $\rk(M) \geq 2$ and let $K \cdot v$ be a non-semisimple symmetric isotropy orbit, where $0 \neq v \in T_oM$. Then there exists $k \in K$ such that $d_ok(\nu_v(K\cdot v)) = T_v(K\cdot v)$ (as linear subspaces of $T_oM$). In particular, $\dim(\nu_v(K\cdot v)) = \dim(T_v(K\cdot v)) = \frac{1}{2}\dim(M)$. 
\end {prop}  
  
\begin {proof} 
For dual symmetric spaces the isotropy representations coincide. Hence we may assume that $M$ is of non-compact type. From Lemma 4.1 in \cite{BO18} we know that $T_v(K\cdot v)$ and $\nu_v(K\cdot v)$ are  complementary Lie triple systems in $T_oM$ and the abelian part of $\nu_v(K\cdot v)$ coincides with $\RR v$. Let $\Sigma $ be the connected, complete, totally geodesic submanifold of $M$ with $T_o\Sigma = \nu_v(K\cdot v)$.  Let $K'$ be the connected subgroup of $K$ with Lie algebra $[\nu_v(K\cdot v),\nu_v(K\cdot v)]$. Since $\RR v$ is the abelian part of $\nu_v(K\cdot v)$, we have $K' \cdot v = \{v\}$ and therefore $K' \subseteq K_v$, where $K_v$ is the isotropy group of $K$ at $v$. Moreover, since $K'$ is connected, we also have $K' \subseteq (K_v)^o$. Let $0 \neq w \in T_v(K\cdot v)$ be tangent to the local Euclidean factor of $K\cdot v$ at $v$. Then $(K_v)^o$ fixes $w$, since $(K,K_v)$ is a symmetric pair, and hence $K'$ fixes $w$ as well. This means that $w$ is a fixed vector of the image of the slice representation of $K'$ on the normal space $\nu_o\Sigma =  T_v(K\cdot v)$. From Lemma 3.2(ii) in \cite {BO16} we see that the abelian part of $T_v (K\cdot v)$ contains $\RR w$. Thus the Lie triple system $T_v(K\cdot v)$ in non-semisimple and so, by Proposition 4.2 in \cite{BO18}, 
\begin{equation}\label{tano}
T_v(K\cdot v) = \nu _w (K\cdot w) 
\end{equation}
and $\RR w$ is the abelian part of $T_v(K\cdot v)$.  Since $w$ is an arbitrary non-zero vector tangent to the local Euclidean factor of $K \cdot v$ at $v$, we conclude that the local Euclidean factor of $K\cdot v$ is $1$-dimensional. 
  
Since $K' \subset (K_v)^o$ and $(K_v)^o$ leaves invariant the Lie triple system $\nu_v(K\cdot v)$, the set of fixed vectors of $(K_v)^o$ in $\nu_v(K \cdot v)$ is $\RR v$, the abelian part of this Lie triple system. The set of fixed vectors of $(K_v)^o$ in $T_v(K\cdot v)$ is $\RR w$, the tangent space to the local Euclidean factor of this symmetric isotropy orbit. Altogether we see that $\VV = \RR v \oplus \RR w$ is the set of fixed vectors of $(K_v)^o$ in $T_oM$.  Clearly, $\VV$ is a $2$-dimensional Lie triple system. However, $\VV$ is not abelian because $w$ is not in the centralizer $Z_\cp(v) = \nu_v(K\cdot v)$ of $v$ in $\cp \cong T_oM$. It follows that $\VV$ is the tangent space at $o$ of a totally geodesic real hyperbolic plane $\RR H^2$ in $M$. Let $\bar{K} \cong SO_2$ be the connected subgroup of $K$ with Lie algebra $\bar{\ck} = [\VV,\VV] \cong \cs\co_2$. We can assume that $\|v\| = \|w\|$. Since $\bar{K}$ acts transitively on the spheres in $T_o\RR H^2$, there exists $k \in \bar{K}$ with $d_ok(v) = w$. Using (\ref{tano}), this implies 
\[
d_ok(\nu_v(K\cdot v)) = d_ok(Z_\cp(v)) = Z_\cp(d_ok(v)) = Z_\cp(w) = \nu _w (K\cdot w)  = T_v(K \cdot v),
\]
which finishes the proof.
\end{proof}

Assume that $K \cdot v$ is a symmetric isotropy orbit in $\cp \cong T_oM$. It was shown in \cite{BO18}, Lemma 4.1, that $T_v(K\cdot v)$ and $\nu _v(K\cdot v)$ are complementary Lie triple systems and the abelian part of $\nu _v(K\cdot v)$ coincides with $\RR v$. We denote by $\Sigma$ and $\Sigma^\perp$ the connected, complete, totally geodesic submanifolds of $M$ with $T_o\Sigma = T_v(K\cdot v)$ and $T_o\Sigma^\perp = \nu_v(K \cdot v)$ respectively. We write $\Sigma = G'/K'$, where $G'$ is the subgroup of $G$ consisting of the glide transformations of $\Sigma$ and $K' = (G')_o$.

\begin {lm} \label{non-semi}
The identity components of $K'$ and $K_v$ coincide, that is, $(K')^o = (K_v)^o$. Consequently, if $M$ is of compact type,  the symmetric spaces  $\Sigma$ and $K \cdot v$ are locally equivalent up to homothety in any local de Rham factor.
\end {lm}

\begin {proof}
Let $G^\Sigma = \{ g\in I(M) : g(\Sigma)  = \Sigma\}$ and $K^\Sigma$ be the isotropy group of $G^\Sigma$ at $o$. Then $K^\Sigma$ leaves $T_o\Sigma^\perp$ invariant and hence its identity component $(K^\Sigma)^o$ fixes $v$ since the abelian part of $\nu _v(K\cdot v)$ coincides with $\RR v$. Thus we have $(K^\Sigma)^o \subseteq (K_v)^o$. On the other hand, $K_v$ leaves $T_o\Sigma = T_v (K \cdot v)$ invariant, which implies $K_v \subseteq K^\Sigma$ and, in particular, $(K_v)^o \subseteq (K^\Sigma)^o$. Altogether this gives $(K_v)^o = (K^\Sigma)^o$.

Assume that $(K')^o$ is a proper subgroup of $(K^\Sigma)^o$. Since $K' $ is an ideal of $K^\Sigma$ and $K^\Sigma$ is compact, there exists a non-trivial connected normal subgroup $H$ of $(K^\Sigma)^o$ acting trivially on $T_o\Sigma = T_v(K\cdot v)$. Then $H$ is a non-trivial subgroup of $K_v$ acting trivially on $T_v(K\cdot v)$. Therefore $H$ acts trivially on the isotropy orbit $K \cdot v$ in $T_oM$ and also in its affine span, say $v + \VV$. Since $K$ acts linearly on $T_oM$, we have $K \cdot \VV \subseteq \VV$, and therefore $\VV = T_oM$ since $K$ acts irreducibly on $T_oM$. So $H$ must be trivial, which is a contradiction, and we conclude that  
$(K')^o = (K_v)^o$.
\end {proof}
   
From Proposition \ref{half} and Lemma \ref{non-semi} we obtain: 
   
\begin {cor} \label{non-semiCOR}
Let $M= G/K$ be an irreducible, simply connected, Riemannian symmetric space and $K \cdot v$ be a symmetric isotropy orbit. Then the following statements are equivalent:
\begin{itemize}
\item[(i)] $K \cdot v$ is a non-semisimple symmetric isotropy orbit in $T_oM$.
\item[(ii)]  $T_v(K \cdot v)$ is a non-semisimple Lie triple system in $\cp$. 
\item[(iii)] There exists $k \in K$ such that $k(\nu_v(K\cdot v)) = T_v(K \cdot v)$.
\end{itemize}
\end {cor}

\begin{proof}
(i) implies (iii) by Proposition \ref{half}. 

Since $K \cdot v$ is a symmetric isotropy orbit, the normal space $\nu_v(K\cdot v)$ is a non-semisimple Lie triple system. If (iii) holds, then $T_v(K \cdot v) = k(\nu_v(K\cdot v))$ is also a non-semisimple Lie triple system. Thus (ii) holds.

Finally, assume that (ii) holds. The tangent space $T_v(K \cdot v)$ and the normal space $\nu_v(K\cdot v)$ of the symmetric isotropy orbit $K \cdot v$ are complementary Lie triple systems. By Lemma \ref{non-semi} we have $(K')^o = (K_v)^o$. Then the symmetric isotropy orbit $K \cdot v$ is non-semisimple since $(K_v)^o = (K')^o$ fixes a non-zero vector of the Lie triple system $T_v(K \cdot v)$, since by assumption the Lie triple system $T_v(K \cdot v)$ is non-semisimple. 
\end{proof}
   
Let $K\cdot v$ be a symmetric isotropy orbit in $T_oM$. Then the tangent space $T_v(K \cdot v)$ and the normal space $\nu _v(K\cdot v)$ are complementary reflective Lie triple systems in $\cp \cong T_oM$. By Theorem 1.2 in \cite {BO16}, $\nu_v(K \cdot v)$ is a maximal Lie triple system in $T_oM$ (the proof of the maximality involves delicate geometric arguments). We will now prove that $T_v(K\cdot v)$ is also a maximal Lie triple system in $T_oM$.
   
\begin {prop} \label{tanMaximal} 
Let $M=G/K$ be an irreducible, simply connected, Riemannian symmetric space with $\rk(M) \geq 2$ and let $K \cdot v\subset T_oM$ be a symmetric isotropy orbit. Then $T_v (K \cdot v)$ is a maximal Lie triple system in $T_oM$.  
\end{prop}
 
\begin {proof}
We can assume that $M$ is of non-compact type. Let $\Sigma = G'/K'$ be the connected, complete, totally geodesic submanifold of $M$ with $T_o\Sigma = T_v(K\cdot v)$, where $G'$ is the subgroup of $G$ consisting of the glide transformations of $\Sigma$. For the corresponding Lie algebras we have $\cg' = [T_o\Sigma,T_o\Sigma] \oplus T_o\Sigma$ and $\ck' = [T_o\Sigma,T_o\Sigma]$. Let $\Sigma^\perp = G''/K''$ be the connected, complete, totally geodesic submanifold of $M$ with $T_o\Sigma^\perp = \nu_v(K\cdot v)$, where $G''$ is the subgroup of $G$ consisting of the glide transformations of $\Sigma^\perp$. 

Now consider the slice representation $\rho : K' \to SO(\nu _v(K\cdot v))$ and the isotropy representation $\chi : K'' \to SO(\nu_v(K\cdot v))$. By Lemma 3.2(i) in \cite{BO16}, $\rho(K')$ is a normal subgroup of $\chi(K'')$. Since $M$ is of non-compact type, $\Sigma$ is simply connected and therefore $K'$ is connected. From Lemma \ref{non-semi} we then have $K' = (K_v)^o$. Moreover, since $\RR v$ is the abelian part of $\nu_v(K \cdot v)$, we also have $(K_v)^o = (K^{\Sigma^\perp})^o$, and thus $K' = (K^{\Sigma^\perp})^o$. The image under the isotropy representation of $(K^{\Sigma^\perp})^o$ on $T_o\Sigma ^\nu = \nu _v(K\cdot v)$ coincides with that of $K''$, and therefore $\rho(K') = \chi (K'')$. 
  
Assume that $\hat\Sigma$ is a proper totally geodesic submanifold of $M$ that properly contains $\Sigma$. Then we can write $T_o\hat\Sigma = T_o\Sigma \oplus \VV$ with $\{0\} \neq \VV \subset \nu _v(K\cdot v)$. Note that $\ck' = [T_o\Sigma,T_o\Sigma] \subseteq [T_o\hat\Sigma,T_o\hat\Sigma]$ and therefore $K'$ leaves $\VV$ invariant. Since $\rho(K') = \chi (K'')$, also $K''$ leaves $\VV$ invariant. This implies that $\Sigma^\perp$ is a Riemannian product $\Sigma^\perp = \Sigma_1 \times \Sigma_2$ with $T_o\Sigma_1 = \VV$. Note that $\Sigma _2$ and $\hat\Sigma$ are complementary reflective totally geodesic submanifolds and that the isotropy group $K_2$ at $o$ of the glide transformations of $\Sigma_2$ acts trivially on $\VV = T_o\Sigma _1$. Thus the image of the slice representation of $K_2$ fixes any element of $\VV$. It follows from Lemma 3.2(ii) in \cite{BO16} that $\Sigma _1$ is flat. 
Since the abelian part of a reflective Lie triple system, in this case $T_o\Sigma^\perp$, is $1$-dimensional, this implies $T_o\Sigma _1 = \RR v$. Hence $\Sigma _2$ is the semisimple part of the reflective submanifold $\Sigma^\perp$. Corollary 3.3 in \cite {BO16} tells us that $\Sigma _2 $ is not reflective, which is a contradiction. It follows that $\Sigma$ is maximal and hence $T_v(K\cdot v)$ is a maximal Lie triple system in $T_oM$.
\end {proof}

\begin{lm} \label{slicesphere}
Let $M = G/K$ be a Riemannian symmetric space of rank $1$. Assume that there exists a totally geodesic submanifold 
$\Sigma = G'/K'$ with $\dim(\Sigma) \geq 2$ such that the connected slice representation $\rho : (K')^o\to SO (\nu _o \Sigma)$ is trivial. Then $M$ has constant curvature.
\end{lm} 

\begin{proof} From the assumption we obtain that $\nu_o\Sigma$ is a reflective Lie triple system, which implies that $\Sigma$ is a reflective submanifold. Assume that $M$ has non-constant curvature. Using duality, we can assume that $M$ is a hyperbolic space over $\CC$, $\HH$ or $\OO$. 

If $M$ is a complex hyperbolic space $\CC H^n = SU_{1,n}/S(U_1U_n)$, then $\Sigma$ is a real hyperbolic space $\RR H^n = SO^o_{1,n}/SO_n$ or a complex hyperbolic space $\CC H^k = SU_{1,k}/S(U_1U_k)$ for some $k \in \{0,\ldots,n-1\}$. In the first case the slice representation of the isotropy group is equivalent to the standard representation of $SO_n$ on $\RR^n$, in the second case it is equivalent to the representation of $S(U_1U_k) \cong U_k$ on $\CC^{n-k}$, where $U_1$ acts canonically and $U_k$ acts trivially. 

If $M$ is a quaternionic hyperbolic space $\HH H^n = Sp_{1,n}/Sp_1Sp_n$, then $\Sigma$ is a complex hyperbolic space $\CC H^n = SU_{1,n}/S(U_1U_n)$ or a quaternionic hyperbolic space $\HH H^k = Sp_{1,k}/Sp_1Sp_k$ for some $k \in \{0,\ldots,n-1\}$. In the first case the slice representation of the isotropy group is equivalent to the standard representation of $S(U_1U_n) \cong U_n$ on $\CC^n$, in the second case it is equivalent to the representation of $Sp_1 Sp_k$ on $\HH^{n-k}$, where $Sp_1$ acts canonically and $Sp_k$ acts trivially.

If $M$ is a Cayley hyperbolic plane $\OO H^2 = F_4^{-20}/Spin_9$, then $\Sigma$ is a quaternionic hyperbolic plane $\HH H^2 = Sp_{1,2}/Sp_1Sp_2$ or a Cayley hyperbolic line $\OO H^1 = Spin_{1,8}/Spin_8 \cong \RR H^8$. In the first case the slice representation of the isotropy group is equivalent to the standard representation of $Sp_1Sp_2$ on $\HH^2$, in the second case it is equivalent to one of the two inequivalent spin representations of $Spin_8$ on $\RR^8$

None of these slice representations is trivial and it follows that $M$ has constant curvature.
\end{proof}

From the Slice Lemma 3.1 in \cite{BO16} and Lemma \ref{slicesphere} we obtain the following interesting result that generalizes Iwahori's result in \cite{Iw66}.

\begin{cor}\label{slicecoro11} Let $M$ be an irreducible Riemannian symmetric space and assume that there exists a non-flat proper totally geodesic submanifold $\Sigma$ of $M$ with trivial connected slice representation. Then $M$ has constant curvature.
\end{cor}

\begin{prop} \label{slicelemmaG} 
Let $M= G/K$ be a simply connected Riemannian symmetric space of non-positive curvature and $\Sigma =G'/K'$ be a semisimple totally geodesic submanifold of $M$ with $o \in \Sigma$, where $G' \subseteq G$ are the glide transformations of $\Sigma$. Let $\Sigma = \Sigma _1 \times \ldots \times \Sigma _l$ and $M = M_0 \times \ldots \times M_g$ be the de Rham decompositions of $\Sigma$ and $M$ respectively, where the Euclidean factor $M_0$ may be trivial. Assume that the slice representation $\rho: K' \to SO(\nu_o\Sigma)$ is trivial. Then $l \leq g$ and, up to a permutation of the indices, $\Sigma _i \subseteq M_i$ for $1 \leq i \leq l$. Moreover, if $\Sigma _i$ is strictly contained in $M_i$, then $M_i$ and $\Sigma _i$ have constant curvature.
\end{prop}

\begin{proof} Let $i\geq 1$ and consider the irreducible de Rham factor $M_i$ of $M$. Let $0 \neq z \in T_oM_i$ and write $z=v+w$ with $v \in T_o\Sigma$ and $w\in \nu_o\Sigma$. If $k'\in K'$, then $k'v +w = k'v +k'w = k'z \in T_oM_i$ and hence $k'z -z= k'v -v \in T_o\Sigma$. The space
\[
v + \mbox{span} \{ k'v- v : k'\in K'\}
\]
coincides with the affine subspace of $T_o\Sigma$ generated by the orbit $K' \cdot v$. This affine subspace must contain $0 \in T_o\Sigma$, because otherwise, the vector $0 \neq u \in v +\mbox{span}  \{ k'v- v:k'\in K'\}$ of minimal distance to $0$ would be fixed by $K'$, which is a contradiction since $\Sigma$ is semisimple and so $K'$ has no fixed non-zero vectors. Then $-v \in \mbox{span}\{ k'v- v:k'\in K'\}$ and hence $v \in \mbox{span} \{ k'v- v:k'\in K'\}\subseteq T_o\Sigma$. But $k'v-v = k'z-z' \in T_oM_i$ for all $k'\in K'$, and therefore $v\in T_oM_i$. This shows that 
\[
T_oM_i = (T_oM_i\cap T_o\Sigma) \oplus  (T_oM_i\cap \nu_o\Sigma) .
\] 
Assume that $T_oM_i \cap T_o\Sigma \neq \{0\}$. A priori, $T_oM_i \cap T_o\Sigma$ is not necessarily irreducible. However, since $T_oM_i \cap T_o\Sigma$ is $K'$-invariant, it is the tangent space of a Riemannian factor $\Sigma_i'$ of $\Sigma$. By assumption, the slice representation of $\Sigma'_i$, considered as a totally geodesic submanifold of the de Rham factor $M_i$ of $M$, is trivial. Observe that $M_i$ is not the Euclidean factor, because otherwise $\Sigma'_i$ would be flat and so $\Sigma$ would be non-semisimple. It follows that $\Sigma'_i$ is irreducible and thus a de Rham factor $\Sigma_i$. It then follows from Corollary \ref{slicecoro11} that $M_i$, and hence also $\Sigma_i$, has constant curvature. 
\end{proof}

If in Proposition \ref{slicelemmaG} the submanifold $\Sigma$ is not semisimple, then the analogous conclusion holds by adding a flat in the factors of $M$ which are not factors of $\Sigma$.

\section {Extrinsic isometries of maximal totally geodesic submanifolds} \label{extiso}
  
Let $\Sigma =G'/K'$ be a connected, complete, totally geodesic submanifold of $M$ with $o \in \Sigma$, where $G' \subset G$ is the group of glide transformations of $\Sigma$ and $K' = (G')_o$. The subgroup $G^\Sigma = \{g \in I(M) : g(\Sigma) = \Sigma\}$ of $I(M)$ is in general neither connected nor effective on $\Sigma$. We always have $\sigma_o \in G^\Sigma$, where $\sigma_o$ is the geodesic symmetry of $M$ at $o$. Note that $G'$ is a normal subgroup of $G^\Sigma$ and that $K'$ is a normal subgroup of $K^\Sigma = (G^\Sigma)_o$.
  
Without loss of generality we may assume that $M$ is of non-compact type. Then $\Sigma$ is simply connected and hence $K'$ must be connected. 

\begin{re} \label{projKilling}
\rm Let $X\in \cg$ and $X^*$ be the corresponding Killing vector field on $M = G/K$ (see Section \ref{pre}). For $p \in \Sigma$ we denote by $X^\Sigma_p$ the orthogonal projection of $X^*_p$ onto $T_p\Sigma$. Then $X^\Sigma$ is a Killing vector field on $\Sigma$. It is well known, and standard to show, that there exists $Z\in \cg' = [T_o\Sigma, T_o\Sigma] \oplus T_o\Sigma \subseteq \cg$ such that $Z^\Sigma = (Z^*)_{|\Sigma} =  X^\Sigma$. The key fact of the argument is to show that a Killing vector field induced by $G$ projects constantly along any flat totally geodesic submanifold.
\end {re}
  
\begin{lm} \label{<}
If $\Sigma$ is maximal and $\dim(\ck') < \dim(\ck^\Sigma)$ (or equivalently, $\dim(\cg') < \dim(\cg^\Sigma)$), then $\Sigma$ is a reflective submanifold of $M$. 
\end{lm}
  
\begin {proof} Let $0 \neq X \in \ck^\Sigma \setminus \ck^\prime$. By Remark \ref{projKilling}, there exists $Z \in \ck'$ such that $Z^\Sigma = X^\Sigma$. By adding $-Z$, we may assume that $X ^\Sigma =0 = {X^*}_{\vert \Sigma}$, where the last equality holds because $X\in \ck'$ and so the restriction of $X^*$ to $\Sigma$ is always tangent to $\Sigma$. 
  
Then $\Exp(tX)$ is a non-trival one-parameter group of isometries of $M$ acting trivially on $\Sigma$ (in particular, fixes $o$ and leaves $\Sigma$ invariant) and $\phi^t = d_o\Exp(tX)$ is a one-parameter group of linear isometries of $T_oM$ with ${\phi^t}_{\vert T_o\Sigma} = \id_{T_o\Sigma}$ for all $t \in \RR$. Since $\Sigma$ is a maximal totally geodesic submanifold of $M$, we have $T_o\Sigma = \{u \in T_oM : \phi^t(u) = u \mbox{ for all } t \in \RR \}$. It follows that the dimension of $\nu _o\Sigma$ is even, say equal to $2d$ for some $0 < d \in \ZZ$. We can find an orthonormal basis $e_1,f_1,\ldots, e_d, f_d$ of $\nu _o\Sigma$ and $a_1\geq a_2 \geq \cdots \geq a_d >0$ such that   
\[
\phi ^t(e_i) = \cos(2\pi a_i t) e_i + \sin(2\pi a_i t) f_i \ \mbox{ and }\ 
\phi ^t(f_i) = -\sin(2\pi a_i t) e_i + \cos(2\pi a_i t)f_i
\]
for all $i \in \{1,\ldots,d\}$.

Assume that $a_1 > a_d$ and put $t_0 = \frac {1}{a_1}$. Then $\VV = \{u \in T_oM : \phi^{t_o}(u) = u\}$ contains $T_o\Sigma \oplus \RR e_1 \oplus \RR f_1$ and is perpendicular to $\RR e_d \oplus \RR f_d$. Let $\Sigma'$ be the connected, complete, totally geodesic submanifold of $M$ with $T_o\Sigma' = \VV$. Then $\Sigma'$ is a proper submanifold of $M$ and properly contains $\Sigma$, which is a contradiction to the maximality of $\Sigma$. It follows that there exists $0 < a \in \RR$ so that $a_1 = \ldots = a_d = a$. Then $\phi^{\frac{1}{2a}}$ is the orthogonal reflection of $T_oM$ in the subspace $T_o\Sigma$. The corresponding isometry $\Exp(\frac {1}{2a}X)$ is the geodesic reflection of $M$ in $\Sigma$. Consequently, $\Sigma$ is a reflective submanifold of $M$.
\end {proof}
  
We will now proceed with the proof of Theorem \ref{main}. The equivalence of (ii) and (iii) is an immediate consequence of Proposition 4.2 in \cite{BO18}. If (iii) holds, then $\Sigma$ is maximal by Proposition \ref{tanMaximal}. Moreover, from Proposition 4.2 in \cite{BO18}, the abelian part of the Lie triple system $\nu_v(K\cdot v) = Z_\cp(v)$ is $1$-dimensional and so it must coincide with $\RR v$. Then $\rho((K')^o)$ must leave $\RR v$ invariant. Since $\rho((K')^o) \subset SO(\nu_v(K\cdot v))$, the group $\rho((K')^o)$ must fix $v$. This proves (i). Note that all these implications do not use the assumption that $\dim(\Sigma) \geq \frac{1}{2} \dim(M)$. 

It remains to prove that (i) implies (ii). We only need to show that $\Sigma$ is reflective. In fact, if $\rho((K')^o)$ fixes $v$, then the complementary reflective submanifold is non-semisimple according to Lemma 3.2(ii) in \cite{BO16}. The proof that $\Sigma$ is reflective if (i) is satisfied will be done in several steps. 

\medskip
{\sc Case 1:} $\Sigma$ is non-semisimple.

\medskip
Since $\Sigma$ is maximal, it follows from Theorem 1.2 in \cite{BO16} that $T_o\Sigma$ is the normal space of a symmetric isotropy orbit. Then, by Lemma 4.1 of \cite{BO18}, $\Sigma$ is reflective. Moreover, from Lemma 3.2(ii) in \cite {BO16} we see that $\nu_o\Sigma$ is non-semisimple and so, by Corollary \ref{non-semiCOR}, $k(\nu_o\Sigma)  = T_o\Sigma$ for some $k\in K$.

\medskip
{\sc Case 2:} $\Sigma$ is simple.

\medskip
By duality, we can assume that $M$ is of compact type. We will use again the Riemannian symmetric space $\bar{M} = G/G_{\pi(o)} = M/{\sim}$ that we encountered in Section \ref{pre}, where $p \sim q$ if $G_p = G_q$ and $\pi : M \to \bar{M}$ is the canonical projection. We put $\bar{o} = \pi(o)$ and identify $T_oM$ with $T_{\bar{o}}\bar{M}$ by means of the isomorphism $d_o\pi : T_oM \to T_{\bar{o}}\bar{M}$. Then we have $(G_{\bar{o}})^o = K$. In this way, the identity component of the isotropy group at $\bar{o}$ of the glide transformations of $\bar{\Sigma} = \pi(\Sigma)$ is canonically identified with the identity component of $K'$. Note that $\bar{\Sigma}$ is a maximal simple totally geodesic submanifold of $\bar{M}$ since $\Sigma$ is a maximal simple totally geodesic submanifold of $M$.

We define $\VV = \{ \xi \in T_{\bar{o}}\bar{M} : \xi \mbox{ is fixed by } K'\}$. From our assumption (i) we have $\VV \neq \{0\}$. Let $\bar{\Sigma}'$ be the connected, complete, compact, totally geodesic submanifold of $\bar{M}$ with $T_{\bar{o}}\bar{\Sigma}' = \VV$. Since $\bar\Sigma$ is simple, $(K')^o$ acts irreducibly on $T_{\bar{o}}\bar\Sigma$, which implies that $\VV$ is perpendicular to $T_{\bar{o}}\bar\Sigma$. 

Since $\bar{\Sigma}'$ is a compact Riemannian symmetric space, there exists a non-trivial closed geodesic $\gamma_{\bar v}(t)$ of (minimal) period $1$ for some $\bar{v} \in \VV$. Let  $g^{\bar v} \in I(\bar{M})$ be as in Corollary \ref{order2}. Since $\bar{v}$ is fixed by $(K')^o$, $g^{\bar v}$ commutes with $(K')^o$ and hence $d_{\bar{o}}g^{\bar{v}}(T_{\bar{o}}\bar{\Sigma})$ is a $(K')^o$-invariant subspace on which $(K')^o$ acts irreducibly. Then $d_{\bar{o}}g^{\bar v}(T_{\bar{o}}\bar{\Sigma})$ is, as well as $T_{\bar{o}}\bar{\Sigma}$, perpendicular to $\bar{v}$. Since $\dim(d_{\bar{o}}g^{\bar{v}}(T_{\bar{o}}\bar{\Sigma})) = \dim(T_{\bar{o}}\bar{\Sigma}) \geq \frac{1}{2} \dim(T_{\bar{o}}\bar{M})$ by assumption, these two subspaces intersect in a non-trivial $(K')^o$-invariant subspace. Since $(K')^o$ acts irreducibly on $T_{\bar{o}}\bar{\Sigma}$, this implies $d_{\bar{o}}g^{\bar{v}}(T_{\bar{o}}\bar{\Sigma}) = T_{\bar{o}}\bar{\Sigma} $, or equivalently, $g^{\bar{v}}(\bar{\Sigma}) = \bar{\Sigma}$. 

Let $E_{\pm 1} \subset T_{\bar{o}}\bar{M}$ be the eigenspaces corresponding to the eigenvalues $\pm 1$ of $d_{\bar{o}}g^{\bar{v}}$. Note that $E_1$ (resp.\ $E_{-1}$) is the set of fixed vectors of $d_{\bar{o}}g^{\bar{v}}$ (resp.\ of $d_{\bar{o}}(\sigma_{\bar{o}} \circ  g^{\bar{v}})$). Since $g^{\bar{v}}$ commutes with $(K')^o$, both eigenspaces are $(K')^o$-invariant. Let $\Sigma_{\pm 1}$ be the connected, complete, totally geodesic submanifold of $\bar{M}$ with $T_{\bar{o}}E_{\pm 1} = E_{\pm 1}$.  By construction $\Sigma_1$ and $\Sigma_{-1}$ are reflective submanifolds of $\bar M$ with $T_{\bar{o}}\Sigma_{\pm 1} = \nu_{\bar{o}}\Sigma_{\mp 1}$. We have $T_{\bar{o}}\bar{M} = E_{1} \oplus  E_{-1}$ and $T_{\bar{o}}\bar{\Sigma} = \bar{E}_{1} \oplus \bar{E}_{-1}$, where $\bar{E}_{\pm 1} = E_{\pm 1} \cap T_{\bar{o}}\bar{\Sigma}$. Since $(K')^o$ acts irreducibly on $T_{\bar{o}}\bar{\Sigma}$ we have either $\bar{E}_1 = \{0\}$ or  $\bar{E}_{-1} = \{0\}$. Assume that $\bar{E}_{-1} = \{0\}$. Then $\bar{E}_1 = T_{\bar{o}}\bar{\Sigma}$ is a proper subset of $E_1$, since $\bar{v} \in E_1$ and $\bar{v}$ is perpendicular to $T_{\bar{o}}\bar{\Sigma}$. Thus $\bar{\Sigma}$ is a proper totally geodesic submanifold of $\Sigma_1$, which contradicts the maximality of $\bar{\Sigma}$. It follows that $\bar{E}_1 = \{0\}$ and thus $T_{\bar{o}}\bar{\Sigma} = \bar{E}_{-1}$.  Since $\bar{\Sigma}$ is maximal, we get $\bar{\Sigma} = \Sigma_{-1}$. This shows that $\bar{\Sigma}$, and hence also $\Sigma$, are reflective submanifolds. 

\medskip
{\sc Case 3:} $\Sigma $ is semisimple but not simple.

\medskip
By duality, we can assume that $M$ is of non-compact type.  Then $\Sigma = \Sigma _1 \times \ldots \times \Sigma_a$, where $a \geq 2$ and $\Sigma _i$ is a simple totally geodesic submanifold of $M$. We put $d_i =\dim(\Sigma_i)$ and arrange the factors so that $d_1 \leq \ldots \leq d_a$. Let $G^i \subset G$ be the group of glide transformations of $\Sigma _i$ and $K^i = (G^i)_o$. Then the group $G'$ of glide transformation of $\Sigma$ is $G' = G^1 \times \ldots \times G^a$ and we have $K' = (G')_o= K^1 \times \ldots \times K^a$. We put $\VV_i = T_o\Sigma _i$. By assumption, there exists $0 \neq v \in \nu_o\Sigma$ which is fixed under the slice representation of $K'$. Let $R$ be the Riemannian curvature tensor of $M$ at $o$.  
  
\begin{lm} \label{inside} 
For all $0 \neq u \in T_o\Sigma$ we have $R_{u,v} \notin \ck'$. 
\end {lm}

\begin{proof} We prove this by contradiction. Assume that there exists $0 \neq u \in T_o\Sigma$ with $R_{u,v} \in \ck'$. Since $K'$ fixes $v$, we have $[\ck',v] = \{0\}$ and thus $R_{u,v}v = 0$. Since $M$ is a Riemannian symmetric space, this implies $R_{u,v} = 0$ and hence $R_{d_ok(u),v} = R_{d_ok(u), d_ok(v)} =d_ok \circ R_{u,v} \circ (d_ok)^{-1} =0$ for all $k \in K'$. Thus, if $\WW$ is the linear span of the orbit $K' \cdot u$ in $T_o\Sigma$, we have $R_{w,v}=0$ for all $w \in \WW$. There exists a non-empty subset $J$ of $\{1,  \ldots,a\}$ so that the $K'$-invariant subspace $\WW$ can be written as  $\WW = \bigoplus_{j\in J}  \VV_j$. Now consider the centralizer $Z_{T_oM}(\WW) = \{ z \in  \cp : [\WW,z] = \{0\}\} = \{ z \in T_oM: R_{\WW,z} =\{0\}\}$ of $\WW$ in $T_oM \cong \cp$. Since $M$ is irreducible, $Z_{T_oM}(\WW)$ is a proper subset of $T_oM$. Obviously, we have $\RR v \oplus \left( \bigoplus_{j \notin J} \VV_j \right) \subset Z_{T_oM}(\WW)$. It follows that $\WW + Z_{T_oM}(\WW)$ is a proper Lie triple system in $T_oM$ containing $T_o\Sigma$ as a proper subset. This contradicts the maximality of $\Sigma$. 
\end {proof}
   
To finish the proof of Theorem \ref{main}, we need to show that $\Sigma$ is reflective. We will prove this by contradiction. Thus, assume that $\Sigma$ is not reflective. Let $0 \neq x \in \VV_1$ and define $B_x = R_{x,v} \in \ck \subset \cs\co(T_oM)$. Then, by Lemma \ref{inside}, we have $B_x \notin \ck'$. Since $R$ is $K$-invariant and $\hat{K}^1 = K^2 \times \ldots \times K^a$ acts trivially on $\VV_1$, we have $[B_x, \hat{\ck}^1] = \{0\}$.  Define the one-parameter subgroup $h^t_x$ of $SO(T_oM)$ by $h^t_x = \exp(tB_x)$ and denote by $H^t_x$ the corresponding one-parameter group in $K$ given by $h^t_x = d_oH^t_x$. By Lemma \ref{inside} and Lemma \ref{<}, and since $\Sigma$ is maximal, $B_x(T_o\Sigma)$ is not contained in $T_o\Sigma$. Thus we have $h_x^t(T_o\Sigma) \neq T_o\Sigma$ for sufficiently small $t \neq 0$. Equivalently, $\Sigma^t = H_x^t(\Sigma) \neq \Sigma$ for sufficiently small $t \neq 0$.  
 
We now consider the totally geodesic submanifold $\hat{\Sigma}_1 = \Sigma _2 \times \ldots \times \Sigma _a$ of $M$. The group $\hat{G}^1$ of glide transformations of $\hat{\Sigma}_1$ is $\hat{G}^1 = G^2 \times \ldots \times G^a$ and the isotropy group at $o$ is $\hat{K}^1 = K^2 \times \ldots \times K^a$. Since $[B_x, \hat{\ck}^1] = \{0\}$,
the totally geodesic submanifold $\hat{\Sigma}_1^t = H^t_x(\hat{\Sigma}_1)$ has the same glide isotropy group $\hat {K}^1$ at $o$ as $\hat{\Sigma}_1^0 = \hat{\Sigma}_1$. This implies $d_ok(T_o\hat{\Sigma}_1^t) = T_o\hat{\Sigma}_1^t$ for all $t$ and all $k \in \hat{K}^1$. Consequently, $\hat{\VV}_1^t = T_o\hat{\Sigma}_1^t$ is a $\hat{K}^1$-invariant Lie triple system in $T_oM$ for all $t$.
 
\begin{lm} \label{inside2} 
If $t \neq 0$ is sufficiently small, then $\hat{\VV}_1^t \cap T_o\Sigma =\{0\}$. 
\end {lm}
 
\begin{proof} The intersection $\hat{\VV}_1^t \cap T_o\Sigma $ is a $\hat{K}^1$-invariant subspace of $T_o\Sigma$ (recall that $\hat{K}^1$ acts trivially on $\VV_1$). Thus there exist a (possibly empty) subset $J^t$ of $\{2 , \ldots,a\}$ and a (possibly trivial) subspace $\VV ^t$ of $\VV _1$ such that $\hat{\VV}_1^t \cap T_o\Sigma = \VV ^t \oplus \left( \bigoplus_{j \in J^t} \VV _j \right)$. Since $\VV_1$ is perpendicular to $\hat{\VV}_1^0 = \hat{\VV}_1 = T_o\hat {\Sigma}_1$, we have $\VV ^t = \{0\}$ for sufficiently small $t$. 

On the one hand, if the intersection $\hat{\VV}_1^t \cap T_o\Sigma$ is non-trivial, it is the tangent space of a product of factors of $\Sigma$, which implies $T_o\Sigma \subseteq  (\hat{\VV}_1^t \cap T_o\Sigma) +  Z_{T_oM}(\hat{\VV}_1^t \cap T_o\Sigma)$. In fact, we have equality here since $\Sigma$ is maximal.

On the other hand, since $\hat{\VV}_1^t \cap T_o\Sigma $ is $\hat{K}^1$-invariant, it is the tangent space of a product of factors of $\hat{\Sigma}_1^t$ (recall that $\hat{K}^1$ is the isotropy group at $o$ of the glide transformations of $\hat{\Sigma}_1^t$ for all $t$). Since $\hat {\Sigma}_1^t$ is a (not necessarily simple) factor of  $\Sigma^t = H_x^t(\Sigma)$, we conclude that this intersection is the tangent space to a (not necessarily simple) factor of $\Sigma^t$. This implies $T_o\Sigma^t \subseteq (\hat{\VV}_1^t \cap T_o\Sigma) +  Z_{T_oM}(\hat{\VV}_1^t \cap T_o\Sigma)$.

Altogether we now see that the Lie triple system $(\hat{\VV}_1^t \cap T_o\Sigma) +  Z_{T_oM}(\hat{\VV}_1^t \cap T_o\Sigma)$ contains $T_o\Sigma + T_o\Sigma ^t$. Recall that, for sufficiently small $t \neq 0$, the subspace  
$T_o\Sigma^t = h_x^t(T_o\Sigma)$ is different from $T_o\Sigma$ and so $T_o\Sigma$ is a proper subspace of $T_o\Sigma + T_o\Sigma^t$. It follows that $T_o\Sigma$ is a proper subspace of the Lie triple system $(\hat{\VV}_1^t \cap T_o\Sigma) +  Z_{T_oM}(\hat{\VV}_1^t \cap T_o\Sigma)$. Since $\Sigma$ is maximal, this implies $(\hat{\VV}_1^t \cap T_o\Sigma) +  Z_{T_oM}(\hat{\VV}_1^t \cap T_o\Sigma) = T_oM$ and thus $\hat{\VV}_1^t \cap T_o\Sigma = \{0\}$ for sufficiently small $t \neq 0$.
\end {proof}

Let $\pi : T_oM \to \nu _o\Sigma$ be the orthogonal projection from $T_oM$ onto the normal space $\nu_o\Sigma$. Since $\nu_o\Sigma$ is $K'$-invariant and $\hat K^1 \subseteq K'$, $\nu_o\Sigma$  is also $\hat{K}^1$-invariant. This implies that $\pi$ is $\hat{K}^1$-equivariant. From Lemma \ref{inside2} we then obtain that $\pi : T_o\hat{\Sigma}_1^t \to \nu_o\Sigma$ is injective for sufficiently small $t\neq 0$. 
 
The totally geodesic submanifold $\hat{\Sigma}_1^t$ of $M$ is isometric to $\hat{\Sigma}_1$ via $H_x^t$ and we have $H_x^t \hat{K}^1 (H_x^t)^{-1} = \hat{K}^1 = K^2 \times \ldots \times K^a$. Furthermore, $T_o\hat{\Sigma}^t_1 = h_x^t(\VV_2) \oplus \ldots \oplus h_x^t (\VV_a)$ and $K^i$ acts irreducibly on $h_x^t(\VV_i)$ and trivially on $h_x^t (\VV_j)$ for $i,j \in \{2, \ldots , a\}$, $i\neq j$. From this and the fact that $\pi$ is $\hat{K}^1$-equivariant it is not hard to see that we have the orthogonal decomposition $\pi(T_o\hat{\Sigma}_1^t) =  \pi(h_x^t(\VV_2)) \oplus 
 \ldots \oplus \pi(h_x^t (\VV_a))$ and that $\hat{K}^1$ acts irreducibly on $\pi(h_x^t(\VV_i))$ and trivially on $\pi(h_x^t (\VV_j))$ for $i,j \in \{2, \cdots , a\}$, $i\neq j$. Moreover, the irreducible representation of $\hat{K}^1$ on $\pi(h_x^t (\VV_i))$ is equivalent to the representation of $\hat{K}^1$ on $\VV_i = T_o\Sigma_i$ for $i\geq2$. 
 
Since $v \neq 0$ is fixed by $\hat{K}^1$, the subspace $\pi(T_o\hat{\Sigma}_1^t)$ is perpendicular to $\RR v$. Let $\WW \subset \nu_o\Sigma$ be the linear span of all $\hat{K}^1$-invariant subspaces of $\nu_o\Sigma$ on which the representation of $\hat K^1$ is equivalent to the representation of $\hat{K}^1$ on $\VV_i$ for some $i\in \{2, \ldots , a\}$. Then $\WW$ is perpendicular to $\RR v$ and contains the subspace $\pi(T_o\hat{\Sigma}_1^t)$ of dimension $d_2 + \ldots + d_a = \dim(\Sigma) -d_1$. If $\pi(T_o\hat{\Sigma}_1^t)$ is a proper subspace of $\WW$, then $\dim(\WW) \geq (d_2 + \ldots + d_a) + d_i$ for some $i \in \{2, \ldots a\}$. Since $d_1\leq \cdots \leq d_a$, we obtain that $\dim(\WW)\geq \dim(\Sigma)$ and so the codimension of $\Sigma$ is at least $\dim(\Sigma) + 1$, which is a contradiction. Hence we have $\WW =  \pi(T_o\hat{\Sigma}_1 ^t)$ and $\hat{K}^1$ acts on $\WW$ as an $s$-representation, equivalent to the isotropy representation of $\hat{K}^1$ on $T_o\hat{\Sigma}_1$. 
 
Next, $K^1$ commutes with $\hat{K}^1$ and leaves $\nu_o\Sigma$ invariant. From the definition of $\WW$ we therefore see that $K^1$ leaves $\WW$ invariant. Moreover, $K^1_{\vert \WW} = \{{k^1}_{\vert \WW} : k^1\in K^1\}$ lies in the centralizer, and so in the normalizer of $\hat{K}^1_{\vert \WW}$ in $SO(\WW)$. Since $K^1_{\vert \WW}$ acts as an $s$-representation, we obtain $ K^1_{\vert \WW} \subseteq \hat{K}^1_{\vert \WW}$ from Lemma 5.2.2 in \cite{BCO16}.

It is well-known that the dimension of the centre of the isotropy group of an irreducible Riemannian symmetric space is either $0$ or $1$ (and in the latter case the space is Hermitian symmetric). If the irreducible Riemannian symmetric space $\Sigma _1$ is not a real hyperbolic plane $\RR H^2$, then $K^1$ is not abelian and $K^1$ does not act effectively on $\WW$. Interchanging $\Sigma _1$ with $\hat{\Sigma}_1$, one can show with similar arguments that there exists a $K^1$-invariant subspace $\tilde{\VV}_1$ of $\nu _o\Sigma$ on which the representation of $K^1$ is equivalent to the one of $K^1$ on $T_o\Sigma _1 = \VV _1$. This subspace satisfies $\tilde{\VV}_1 \cap (\RR v \oplus \WW) = \{0\}$. Since $\dim(\WW) = \dim(\Sigma) - d_1$, we conclude that $\dim(\nu_o\Sigma) \geq \dim(\Sigma) + 1$, which is a contradiction. We conclude that $\Sigma _1$ is a real hyperbolic plane $\RR H^2$ and therefore $\dim(\WW) = \dim(\Sigma) -2$. 

This implies $\dim(\nu_o\Sigma) \geq \dim(\Sigma) -1$, since $v$ is perpendicular to $\WW$. Thus there exists $w \in \nu _o\Sigma$, possibly $w = 0$ if $\dim(\nu_o\Sigma) = \dim(\Sigma) -1$, perpendicular to $\RR v \oplus \WW$ such that $\nu _o\Sigma = \RR w \oplus \RR v \oplus \WW$. The group $K^1$ fixes $w$. Let $k \in K^1$ be non-trivial. Then there exists $\hat{k} \in \hat{K}^1$ such that $k\hat{k}^{-1}$ acts trivially on $\WW$ and hence on 
$\nu _o\Sigma$. Note that, since  $\hat{k}$ acts trivially on $\Sigma _1$, $k\hat{k}^{-1}$ is a non-trivial element in $G^\Sigma$ and lies in the kernel of the slice representation (recall that the image of an $s$-representation coincides with its own connected normalizer in the full orthogonal group, see e.g. Lemma 5.2.2 in \cite {BCO16}). Then, by Proposition \ref {fullslice}, $\Sigma$ is reflective, which contradicts our assumption that $\Sigma$ is non-reflective. This finishes the proof of Theorem \ref {main}.

\section{The index of exceptional Riemannian symmetric spaces} \label{indxexc}

Let $M=G/K$ be an $n$-dimensional irreducible Riemannian symmetric space of non-compact type and denote by $r$ the rank of $M$. Let $\Sigma = G'/K'$ be a connected totally geodesic submanifold of $M$ with codimension $d \geq 1$ and denote by $r_\Sigma$ the rank of $\Sigma$. We can assume that $o \in \Sigma$ and $G^\prime \subseteq G$ is the group of glide transformations of $\Sigma$. Then we have $\dim(G^\prime) - \dim(K^\prime) = \dim(\Sigma) =  n-d$. The dimension of a principal orbit of the isotropy action of $K'$ on $\Sigma$ is $\dim(K^\prime) - \dim(K'_0)  = (n-d) - r_\Sigma$, where $K'_0$ is the principal isotropy group of $K'$ (that is, the isotropy group of $K'$ at a point in a principal orbit of the $K'$-action on $\Sigma$). Altogether this implies 
\begin {equation} \label{dimGprime}
\dim(G') = 2(n-d) - r_\Sigma + \dim(K'_0) .
\end {equation}

Let $M^\ast = G^\ast/K$ and $\Sigma^\ast = G'^\ast/K'$ be the Riemannian symmetric spaces of compact type that are dual to $M = G/K$ and $\Sigma = G'/K'$ respectively. Let  $i(G^\ast)$ denote the index of $G^\ast$, where the compact Lie group $G^\ast$ is considered as a Riemannian symmetric space of compact type. Using \cite{BO17} we get the inequality
\[
\dim (G) -\dim (G') = \dim (G^\ast) -\dim (G'^\ast) \geq i(G^\ast), 
\]
or equivalently, 
\[
n + \dim(K)  -\dim(G') \geq i(G^\ast).
\] 
Using (\ref{dimGprime}), we obtain 
\[
n + \dim(K)  - 2(n-d) -\dim(K'_0) + r_\Sigma \geq i(G^\ast),
\]
or equivalently,
\begin {equation}  \label{refinement}
d\geq \frac{1}{2}( i (G^\ast)  + n - r_\Sigma -\dim(K) + \dim(K'_0)),
\end {equation}
or equivalently,  
\begin {equation}  \label{refinement2}
d\geq \frac{1}{2}( i (G^\ast)  + n - r -\dim(K) + (r-r_\Sigma) + \dim(K'_0)),
\end {equation}
Since $r -r_\Sigma$ and $\dim(K'_0)$ are non-negative, the previous equation implies
\begin {equation} \label{refinement-bound}
d \geq \frac{1}{2}( i(G^\ast)  + n - r -\dim(K) ).
\end {equation} 

We introduce some notations. By $\ell_\Sigma$ we denote the dimension of the principal isotropy algebra $\ck'_0$ of $\Sigma = G'/K'$, thus $\ell_\Sigma = \dim(K_0^\prime)$. We define the integers $\Omega_M$ and $\Lambda^M_\Sigma$ by
\[
\Omega_M = i(G^\ast) + \dim(M) - \rk(M) -\dim(K) \mbox{ and } \Lambda^M_\Sigma = (\rk(M) - \rk(\Sigma)) + \ell_\Sigma.
\]
We will frequently use some data about symmetric spaces, which we summarize in Table \ref{symspacedata}.

\begin{table}[ht]
\caption{Dimension, rank and $\ell$-number of symmetric spaces $\Sigma$} 
\label{symspacedata} 
{\footnotesize\begin{tabular}{ | p{3.2cm}  p{2cm}  p{1cm} p{2cm}  p{1.9cm}  |}
\hline \rule{0pt}{4mm}
\hspace{-1mm}$\Sigma$ & $\dim(\Sigma)$ & $\rk(\Sigma)$ & $\ell_\Sigma$ & Comments \\[1mm]
\hline \rule{0pt}{4mm}
\hspace{-2mm} 
$SL_{r+1}(\bbr)/SO_{r+1}$ & $\frac{1}{2}r(r+3)$ & $r$ & $0$ & $r \geq 2$ \\
$SL_{r+1}(\bbc)/SU_{r+1}$ & $r(r+2)$ & $r$ & $r$ & $r \geq 2$ \\
$SU^*_{2r+2}/Sp_{r+1}$ & $r(2r+3)$ & $r$ & $3(r+1)$ & $r \geq 2$ \\
$E_6^{-26}/F_4$ & $26$ & $2$ & $28$ &  \\
$SO^o_{1,k+1}/SO_{k+1}$ & $k+1$ & $1$ & $\frac{1}{2}(k-1)k$ & $k \geq 1$ \\[1mm]
\hline \rule{0pt}{4mm}
\hspace{-2mm} 
$SO^o_{r,r+k}/SO_{r}SO_{r+k}$ & $r(r+k)$ & $r$ & $\frac{1}{2}(k-1)k$ & $r \geq 2, k \geq 1$ \\
$SO_{2r+1}(\bbc)/SO_{2r+1}$ & $r(2r+1)$ & $r$ & $r$ & $r \geq 2$ \\[1mm]
\hline \rule{0pt}{4mm}
\hspace{-2mm} 
$Sp_r(\bbr)/U_r$ & $r(r+1)$ & $r$ & $0$ & $r \geq 3$ \\
$SU_{r,r}/S(U_rU_r)$ & $2r^2$ & $r$ & $r-1$ & $r \geq 3$ \\
$Sp_r(\bbc)/Sp_r$ & $r(2r+1)$ & $r$ & $r$ & $r \geq 3$ \\
$SO^*_{4r}/U_{2r}$ & $2r(2r-1)$ & $r$ & $3r$ &  $r \geq 3$ \\
$Sp_{r,r}/Sp_rSp_r$ & $4r^2$ & $r$ & $3r$ & $r \geq 3$ \\
$E_7^{-25}/E_6U_1$ & $54$ & $3$ & $28$ &  \\[1mm]
\hline \rule{0pt}{4mm}
\hspace{-2mm} 
$SO^o_{r,r}/SO_{r}SO_{r}$ & $r^2$ & $r$ & $0$ &  $r \geq 4$ \\
$SO_{2r}(\bbc)/SO_{2r}$ & $r(2r-1)$ & $r$ & $r$ &  $r \geq 4$ \\[1mm]
\hline \rule{0pt}{4mm}
\hspace{-2mm} 
$SU_{r,r+k}/S(U_rU_{r+k})$ & $2r(r+k)$ & $r$ & $k^2 + r - 1$ &  $r \geq 2, k \geq 1$ \\
$Sp_{r,r+k}/Sp_rSp_{r+k}$ & $4r(r+k)$ & $r$ & $k(2k+1) + 3r$ &  $r \geq 2, k \geq 1$  \\
$SO^*_{4r+2}/U_{2r+1}$ & $2r(2r+1)$ & $r$ & $3r+1$ & $r \geq 2$  \\
$E_6^{-14}/Spin_{10}U_1$ & $32$ & $2$ & $16$ &   \\
$F_4^{-20}/Spin_9$ & $16$ & $1$ & $21$ &   \\[1mm]
\hline \rule{0pt}{4mm}
\hspace{-2mm}  
$E_6^6/Sp_4$ &  $42$ & $6$ & $0$ & \\
$E_6(\bbc)/E_6$ & $78$ & $6$  & $6$ & \\[1mm]
\hline \rule{0pt}{4mm}
\hspace{-2mm} 
$E_7^7/SU_8$ & $70$ & $7$ & $0$ & \\
$E_7({\mathbb C})/E_7$ & $133$ & $7$ & $7$ & \\[1mm]
\hline \rule{0pt}{4mm}
\hspace{-2mm} 
$E_8^8/SO_{16}$ & $128$ & $8$ & $0$ & \\
$E_8(\bbc)/E_8$ & $248$ & $8$ & $8$ & \\[1mm]
\hline \rule{0pt}{4mm}
\hspace{-2mm} 
$F_4^4/Sp_3Sp_1$ & $28$ & $4$ & $0$ & \\
$E_6^2/SU_6Sp_1$ & $40$ & $4$ & $2$ &  \\
$F_4(\bbc)/F_4$ & $52$ & $4$ & $4$ &  \\
$E_7^{-5}/SO_{12}Sp_1$ & $64$ & $4$ & $9$ & \\
$E_8^{-24}/E_7Sp_1$ & $112$ & $4$ & $28$ & \\[1mm]
\hline \rule{0pt}{4mm}
\hspace{-2mm} 
$G^2_2/SO_4$ & $8$ & $2$ & $0$ & \\
$G_2(\bbc)/G_2$ & $14$ & $2$ & $2$ & \\[1mm]
\hline
\end{tabular}}
\end{table}

Inequalities (\ref{refinement2}) and (\ref{refinement-bound}) give the following estimates for the codimension of $\Sigma$:

\begin{prop} \label{firstestimate}
Let $M=G/K$ be an irreducible Riemannian symmetric space of non-compact type and $\Sigma$ be a connected totally geodesic submanifold of $M$ with $\codim(\Sigma) \geq 1$. Then
\[
\codim(\Sigma) \geq \frac{1}{2}(  \Omega_M + \Lambda^M_\Sigma  ) \geq \frac{1}{2} \Omega_M.
\]
\end{prop}

Using the second inequality in Proposition \ref{firstestimate} we can confirm the conjecture for some symmetric spaces where it was previously unknown.

\begin{thm} \label{easythm} For the following symmetric spaces the index coincides with the reflective index:
\begin{itemize}
\item[(i)] For $M = E_6^2/SU_6Sp_1$ we have $i(M) = 12$.
\item[(ii)] For $M = E_7^7/SU_8$ we have $i(M) = 27$.
\item[(iii)] For $M = E_8^8/SO_{16}$ we have $i(M) = 56$.
\item[(iv)] For $M = Sp_r(\RR)/U_r$ ($r \geq 3$) we have $i(M) = 2(r-1)$.
\end{itemize}
\end{thm}

\begin{proof}
The index of compact simple Lie groups was calculated in \cite{BO17} and the reflective index of irreducible Riemannian symmetric spaces $M$ was calculated in \cite{BO16}. Using these results we obtain:

If $M = G/K = E_6^2/SU_6Sp_1$, then $i(E_6) = 26$, $\dim(M) = 40$, $\rk(M) = 4$, $\dim(K) = 38$. This gives $\codim(\Sigma) \geq  \frac{1}{2}(26 + 40 - 4 - 38) = 12 = i_r(M)$. 

If $M = G/K = E_7^7/SU_8$, then $i(E_7) = 54$, $\dim(M) = 70$, $\rk(M) = 7$, $\dim(K) = 63$. This gives $\codim(\Sigma) \geq  \frac{1}{2}(54 + 70 - 7 - 63) = 27 = i_r(M)$. 

If $M = G/K = E_8^8/SO_{16}$, then $i(E_8) = 112$, $\dim(M) = 128$, $\rk(M) = 8$, $\dim(K) = 120$. This gives $\codim(\Sigma) \geq  \frac{1}{2}(112 + 128 - 8 - 120) = 56 = i_r(M)$.

If $M = G/K = Sp_r(\RR)/U_r$, then $i(Sp_r) = 4(r-1)$, $\dim(M) = r(r+1)$, $\rk(M) = r$, $\dim(K) = r^2$. This gives $\codim(\Sigma) \geq  \frac{1}{2}(4r-4 + r^2 + r - r - r^2) = 2r-2 = i_r(M)$.
\end{proof}

This finishes the proof of Theorem \ref{mainsp}. To finish the proof of Theorem \ref{mainexc}, we need to investigate further the symmetric spaces $E_6^6/Sp_4$, $E_7^{-5}/SO_{12}Sp_1$, $E_7^{-25}/E_6U_1$ and $E_8^{-24}/E_7Sp_1$. We will investigate these four spaces individually, but first summarize in the following corollary what the first inequality in Proposition \ref{firstestimate} tells us in this situation.

\begin{cor} \label{intermediate}
For the index $i(M)$ of $M$ and we have:
\begin{itemize}
\item[(i)] If $M = E_6^6/Sp_4$, then $13 \leq i(M) \leq i_r(M) = 14$.
\item[(ii)] If $M = E_7^{-5}/SO_{12}Sp_1$, then $23 \leq i(M) \leq i_r(M) = 24$.
\item[(iii)] If $M = E_7^{-25}/E_6U_1$, then $13 \leq i(M) \leq i_r(M) = 22$.
\item[(iv)] If $M = E_8^{-24}/E_7Sp_1$, then $42 \leq i(M) \leq i_r(M) = 48$.
\end{itemize}
Moreover, if $i(M) < i_r(M)$ and $\Sigma$ is a totally geodesic submanifold of $M$ with $\codim(\Sigma) = i(M)$, then
\begin{itemize}
\item[(i)] $\Lambda^M_\Sigma = 0$, if $M = E_6^6/Sp_4$.
\item[(ii)] $0 \leq \Lambda^M_\Sigma \leq 1$, if $M = E_7^{-5}/SO_{12}Sp_1$.
\item[(iii)] $0 \leq \Lambda^M_\Sigma \leq 16$, if $M = E_7^{-25}/E_6U_1$.
\item[(iv)] $0 \leq \Lambda^M_\Sigma \leq 10$, if $M = E_8^{-24}/E_7Sp_1$.
\end{itemize}
\end{cor}

\begin{proof}
If $M = G/K = E_6^6/Sp_4$, then $i(E_6) = 26$, $\dim(M) = 42$, $\rk(M) = 6$, $\dim(K) = 36$. This gives $\codim(\Sigma) \geq  \frac{1}{2}(26 + 42 - 6 - 36) = 13 < 14 = i_r(M)$.

If $M = G/K = E_7^{-5}/SO_{12}Sp_1$, then $i(E_7) = 54$, $\dim(M) = 64$, $\rk(M) = 4$, $\dim(K) = 69$. This gives $\codim(\Sigma) \geq  \frac{1}{2}(54 + 64 - 4 - 69) = 22.5 <  24 = i_r(M)$.

If $M = G/K = E_7^{-25}/E_6U_1$, then $i(E_7) = 54$, $\dim(M) = 54$, $\rk(M) = 3$, $\dim(K) = 79$. This gives $\codim(\Sigma) \geq  \frac{1}{2}(54 + 54 - 3 - 79) = 13 <  22 = i_r(M)$. 

If $M = G/K = E_8^{-24}/E_7Sp_1$, then $i(E_8) = 112$, $\dim(M) = 112$, $\rk(M) = 4$, $\dim(K) = 136$. This gives $\codim(\Sigma) \geq  \frac{1}{2}(112 + 112 - 4 - 136) = 42 <  48 = i_r(M)$.

The statements about $\Lambda^M_\Sigma$ follow immediately from the inequality in Proposition \ref{firstestimate} and the above calculations.
\end{proof}

We now proceed with individual arguments for the four remaining exceptional symmetric spaces.

\begin{thm} \label{E6Sp4}
For $M = E_6^6/Sp_4$ we have $i(M) = i_r(M) = 14$.
\end{thm}

\begin{proof}
We already know from \cite{BO16} that $i_r(M) = 14$ and from Corollary \ref{intermediate} that $i(M) \geq 13$. Let $\Sigma$ be a maximal totally geodesic submanifold of $M$ and assume that $\codim(\Sigma) = 13$. Using Corollary \ref{intermediate} we obtain $\rk(M) - \rk(\Sigma) + \ell_\Sigma = \Lambda^M_\Sigma  = 0$. Consequently, we have $\rk(\Sigma) = \rk(M) = 6$ and $\ell_\Sigma = 0$. The maximal totally geodesic submanifolds of maximal rank in $E_6^6/Sp_4$ were classified by Chen and Nagano in \cite{CN78}, and independently by Ikawa and Tasaki in \cite{IT00}. Up to congruency, there are only three such submanifolds, namely
\begin{itemize}
\item[(i)] $\Sigma = \RR \times SO^o_{5,5}/SO_5SO_5$, which is reflective and $\codim(\Sigma) = 16$;
\item[(ii)] $\Sigma = \RR H^2 \times SL_6(\RR)/SO_6$, which is reflective and $\codim(\Sigma) = 20$;
\item[(iii)] $\Sigma = SL_3(\RR)/SO_3 \times SL_3(\RR)/SO_3 \times SL_3(\RR)/SO_3$, which is non-reflective and $\codim(\Sigma) = 27$.
\end{itemize}
In all three cases we have $\codim(\Sigma) > 14 = i_r(M)$, which is a contradicition. Thus we can conclude that $i(M) = i_r(M) = 14$.
\end{proof}

For the remaining three exceptional symmetric spaces we will first prove a theoretical result that will allow us to reduce the number of cases that need to be considered. The following lemma is a slight generalization of a result by Iwahori \cite{Iw66} which states that if an irreducible Riemannian symmetric space admits a totally geodesic hypersurface, then it must be a space of constant curvature. 

\begin{lm} \label{cod1-semisimple} 
Let $M$ be a Riemannian symmetric space with de Rham decomposition $M = M_0 \times M_1 \times \ldots \times M_g$, where $M_0$ is a {\rm (}possibly $0$-dimensional{\rm )} Euclidean space and $M_1,\ldots,M_g$, $g \geq 1$, are irreducible Riemannian symmetric spaces. Let $\Sigma$ be a totally geodesic hypersurface of $M$. Then $\Sigma = M_0 \times M_1 \times \ldots M_{j-1} \times \Sigma_j \times M_{j+1} \times \ldots \times M_g$, where $\Sigma_j$ is a totally geodesic hypersurface of $M_j$ and $M_j$ is a space of constant curvature for some $j \in \{0,\ldots,g\}$. 
\end{lm}

\begin{proof}
We can assume that all spaces contain the base point $o \in M$. The intersection $T_o\Sigma \cap T_oM_j$ is a Lie triple system in $T_oM_j$ of dimension equal to $\dim(M_j) - 1$ or $\dim(M_j)$. Iwahori's result therefore implies that $T_oM_j \subseteq T_o\Sigma$ if $M_j$ is not of constant curvature, or equivalently, $M_j \subseteq \Sigma$ if $M_j$ is not of constant curvature. We can therefore assume that $M_j$ has constant curvature for all $j \in \{1,\ldots,g\}$. Denote by $R^j$ the Riemannian curvature tensor of $M_j$, by $R$ the Riemannian curvature tensor of $M$, and by $\kappa_j$ the constant sectional curvature of $M_j$. Let $X \in T_o\Sigma \cap T_oM_j$ be a unit vector and $Y = Y_0 + \ldots + Y_g \in T_o\Sigma$ with $Y_j  \in T_oM_j$ and $\langle X,Y \rangle = 0$. Then we have $\kappa_jY_j = R^j(Y_j,X)X = R(Y,X)X \in T_o\Sigma \cap T_oM_j$ and therefore $Y_j \in T_o\Sigma \cap T_oM_j$ for all $j \in \{1,\ldots,g\}$. This implies $T_o\Sigma = (T_o\Sigma \cap T_oM_0) \oplus \ldots \oplus (T_o\Sigma \cap T_oM_g)$ and hence $T_oM_j \subseteq T_o\Sigma$ for all but one index $j \in \{0,\ldots,g\}$. For this $j$ we define $\Sigma_j$ to be the totally geodesic hypersurface of $M_j$ corresponding to the Lie triple system $T_o\Sigma \cap T_oM_j$ in $T_oM_j$. This implies the assertion.
\end{proof}

\begin{prop} \label{hyperbolicfactors} 
Let $\Sigma$ be a reducible maximal totally geodesic submanifold of an irreducible Riemannian symmetric space $M$ of non-compact type. Assume that the de Rham decomposition of $\Sigma$ contains a real hyperbolic space $\RR H^k$ ($k \geq 2$), a complex hyperbolic space $\CC H^k$ ($k \geq 2$), the symmetric space $SL_3(\RR)/SO_3$, or the symmetric space $SO^o_{2,2+k}/SO_2SO_{2+k}$ ($k \geq 1$ odd). Then either $\Sigma = \RR H^{k_1} \times \RR H^{k_2}$ for some $k_1,k_2 \geq 2$, or there exists a reflective submanifold $\Sigma'$ of $M$ with $\dim(\Sigma') \geq \dim(\Sigma)$.
\end {prop}

\begin{proof}
As usual, we write $M = G/K$ with $G = I^o(M)$ and $K = G_o$ with $o \in M$. We can assume that $o \in \Sigma$ and write $\Sigma = G'/K'$, where $G' \subset G$ is the group of glide transformations of $\Sigma$ and $K' = (G')_o$. We denote the de Rham factor specified in the assertion by $N$. We can assume that $o \in N$ and write $N = G''/K''$, where $G'' \subset G$ is the group of glide transformations of $N$ and $K'' = (G'')_o$.

There exists a $2$-dimensional reflective submanifold $P$ of $N$ with $o \in P$ such that the geodesic reflection $\tau$ in its perpendicular reflective submanifold $P^\perp$ at $o$ is inner, i.e. $\tau \in K'' \subseteq K' \subset K$. Explicitly, 
\begin{itemize}
\item[(i)] If $N= \RR H^k$, then $P = \RR H^2$ and $P^\perp = \RR H^{k-2}$; 
\item[(ii)] If $N = \CC H^k$, then $P = \CC H^1$ and $P^\perp = \CC H^{k-1}$;
\item[(iii)] If $N = SL_3(\RR)/SO_3$, then $P = \RR H^2$ and $P^\perp = \RR H^2 \times \RR$;
\item[(iv)] If $N = SO^o_{2,2+k}/SO_2SO_{2+k}$, then $P = \CC H^1$ and\\ $P^\perp = SO^o_{2,2+(k-1)}/SO_2SO_{2+(k-1)}$.
\end{itemize}

If $\Sigma$ is non-semisimple, then $\Sigma$ is reflective by Corollary 4.4 in \cite{BO16} and we can choose $\Sigma' = \Sigma$. We assume from now on that $\Sigma$ is semisimple and write $\Sigma = N \times \bar\Sigma$ with a semisimple totally geodesic submanifold $\bar\Sigma$ of $M$ containing $o$.

Let $T$ be the closure of the subgroup of $K \subseteq SO(T_oM)$ (via the isotropy representation) that is generated by $\tau$. Note that $d_o\tau$ is the identity on $T_oP^\perp \oplus T_o\bar\Sigma$ and minus the identity on $T_oP$. In particular, the cardinality of $\{ g|_{T_o\Sigma} : g \in T\}$ is equal to $2$. 

Assume that $\dim(T) > 0$. Then the kernel $T_\Sigma$ of the Lie group homomorphism $T \to SO(T_o\Sigma),\ g \mapsto g|_{T_o\Sigma}$ has positive dimension. Since $K^\prime$ acts almost effectively on $\Sigma$, we cannot have $T_\Sigma \subseteq K^\prime$. Thus for the corresponding Lie algebras we have $\ct_\Sigma \not\subseteq \ck^\prime$ and $\ct_\Sigma \subseteq \cg^\Sigma$. It then follows from Lemma \ref{<} that $\Sigma$ is reflective and we can choose $\Sigma' = \Sigma$.

From now on we assume that $\dim(T) = 0$, or equivalently, that $\tau$ has finite order. By construction, $\tau$ has even order of the form $2^sq$ with $0 < q \in \ZZ$ odd and $0 < s \in \ZZ$. We replace $\tau$ by $\tau^q$. Then $\tau$ has order $2^s$. If $s > 1$, $\tau^{2^{s-1}}$ is involutive and its set of fixed vectors must coincide with  $T_o\Sigma$ since $\Sigma$ is maximal. This implies that $\Sigma$ is reflective and we can choose $\Sigma' = \Sigma$.

Thus we are left with the case that $\tau$ is an involution. If the set of fixed vectors of $\tau$ in the normal space $\nu_o\Sigma$ is trivial, then $P^\perp \times \bar\Sigma \subset \Sigma$ is reflective and it follows from Corollary \ref{proper-reflective} that $\Sigma$ is reflective and we can choose $\Sigma' = \Sigma$.

So let us assume that the subspace $\VV$ of fixed vectors of $\tau$ in $\nu_o\Sigma$ satisfies $\dim(\VV) > 0$. Since $\tau$ is involutive, the totally geodesic submanifold $\Sigma'$ of $M$ with $o \in \Sigma'$ and $T_o\Sigma' = T_oP^\perp \oplus T_o \bar\Sigma \oplus  \VV$ is reflective. If $\dim(\VV) \geq 2$, then $\dim(\Sigma') \geq \dim(\Sigma)$. 

Thus it remains to analyze the case $\dim (\VV) =1$.  Then $\hat\Sigma = P^\perp \times \bar{\Sigma}$ is a totally geodesic hypersurface of $\Sigma'$. Recall that $P^\perp$ is irreducible unless $N = SL_3(\RR)/SO(3)$, where $P^\perp = \RR H^2 \times \RR$. Let $\hat{\Sigma} = \hat{\Sigma}_0 \times \ldots \times \hat{\Sigma}_g$ be the de Rham decomposition of $\hat\Sigma$, where $\hat\Sigma_0$ is the (possibly $0$-dimensional) Euclidean factor. We can arrange the indices so that $\hat\Sigma_g = P^\perp$ unless $N = SL_3(\RR)/SO(3)$, in which case $\hat\Sigma_g = \RR H^2$ and the $1$-dimensional Euclidean factor in $P^\perp = \RR H^2 \times \RR$ coincides with $\hat\Sigma_0$.
Then, by Lemma \ref{cod1-semisimple}, we have $\Sigma' = \hat\Sigma_0 \times \hat\Sigma_1 \times \ldots \hat\Sigma_{j-1} \times \Sigma'_j \times \hat\Sigma_{j+1} \times \ldots \times \hat\Sigma_g$, where $\hat\Sigma_j$ is a totally geodesic hypersurface of $\Sigma'_j$ and $\Sigma'_j$ is a space of constant curvature for some $j \in \{0,\ldots,g\}$. 

Assume that $g \geq 3$. Then there exists $i \in \{1,\ldots,g-1\}$ with $i \neq j$. The subspace ${\mathcal Z} = Z_{T_oM}(T_o\hat\Sigma_i) \oplus T_o\hat\Sigma_i$ of $T_oM$ is a Lie triple system in $T_oM$ containing both $T_o\Sigma$ and $T_o\Sigma'$. Moreover, ${\mathcal Z}$ is a proper subset of $T_oM$ since $M$ is irreducible. By construction, $T_o\Sigma'$ is not contained in $T_o\Sigma$ and therefore $T_o\Sigma$ is a proper subset of ${\mathcal Z}$, which contradicts the maximality of $\Sigma$. Consequently, we have $g \in \{1,2\}$. If $g=1$, then $\bar\Sigma$ is trivial, which contradicts the assumption that $\Sigma = N \times \bar\Sigma$ is reducible. Consequently, we have $g = 2$
and therefore $\hat{\Sigma} = \hat{\Sigma}_0 \times \hat{\Sigma}_1 \times \hat\Sigma_2$.

If $j \in \{0,2\}$, we consider the subspace $Z_{T_oM}(T_o\hat\Sigma_1) \oplus T_o\hat\Sigma_1$ of $T_oM$, which again is a Lie triple system in $T_oM$ containing both $T_o\Sigma$ and $T_o\Sigma'$. Then the same argument as in the previous paragraph leads to a contradiction. We therefore must have $j =1$. Then $\bar\Sigma = \hat{\Sigma}_1 = \RR H^{k_1}$ and $\Sigma'_1 = \RR H^{k_1+1}$ for some $k_1 \geq 2$, since $\hat{\Sigma}_1$ is a non-flat totally geodesic hypersurface of the irreducible space $\Sigma'_1$ of constant negative curvature. It follows that $\Sigma = N \times \RR H^{k_1}$. Finally, following the same arguments as above with $\RR H^{k_1}$ instead of $N$, we conclude that there exists a reflective submanifold $\Sigma'$ of $M$ with $\dim(\Sigma') \geq \dim(\Sigma)$, or $\Sigma = \RR H^{k_1} \times \RR H^{k_2}$ for some $k_1,k_2 \geq 2$.
\end{proof}

We will use Proposition \ref{hyperbolicfactors} to determine the index of $E_7^{-5}/SO_{12}Sp_1$, $E_7^{-25}/E_6U_1$ and $E_8^{-24}/E_7Sp_1$. We frequently need the values of $\ell_\Sigma$ when $\Sigma$ is a hyperbolic space and therefore list them here:
\begin{align*}
\ell_{\RR H^k} & = \dim(\cs\co_{k-1}) = \frac{1}{2}(k-2)(k-1),\\
\ell_{\CC H^k} & = \dim(\cs\cu_{k-1} \oplus \cu_1) = (k-1)^2,\\
\ell_{\HH H^k} & = \dim(\cs\cp_{k-1} \oplus \cs\cp_1) = (k-1)(2k-1) + 3,\\
\ell_{\OO H^2} & = \dim(\cs\co_7) = 21.
\end{align*}
We will also use frequently the fact that $\ell_\Sigma = \ell_{\Sigma_1} + \ldots + \ell_{\Sigma_g}$ if $\Sigma$ is a Riemannian product $\Sigma = \Sigma_1 \times \ldots \times \Sigma_g$.

\begin{thm} \label{E7SO12Sp1}
For $M = E_7^{-5}/SO_{12}Sp_1$ we have $i(M) = i_r(M) = 24$. 
\end{thm}

\begin{proof}
We already know from \cite{BO16} that $i_r(M) = 24$ and from Corollary \ref{intermediate} that $i(M) \geq 23$. Let $\Sigma$ be a maximal totally geodesic submanifold of $M$ and assume that $\codim(\Sigma) = 23$, that is $\dim(\Sigma) = 41$. Using Corollary \ref{intermediate} we obtain $(\rk(M) - \rk(\Sigma)) + \ell_\Sigma = \Lambda^M_\Sigma  \leq 1$. Since $\rk(M) = 4$, we have $(\rk(\Sigma),\ell_\Sigma) \in \{(4,0),(3,0),(4,1)\}$. Since there exist no irreducible symmetric spaces with rank $3$ or $4$ and of dimension $41$, $\Sigma$ must be reducible. 

Assume that the de Rham decomposition of $\Sigma$ contains a rank one factor $\Sigma_1$. Since $\ell_{\Sigma_1} \leq 1$, we have $\Sigma_1 \in \{\RR H^2,\RR H^3,\CC H^2\}$. It then follows from Proposition \ref{hyperbolicfactors} that there exists a reflective submanifold $\Sigma'$ of $M$ with $\dim(\Sigma') \geq \dim(\Sigma)$, which is a contradiction. 

Finally, assume that $\rk(\Sigma) = 4$ and $\Sigma = \Sigma_1 \times \Sigma_2$ with $\rk(\Sigma_i) = 2$ and $\Sigma_i$ irreducible. Since $\dim(\Sigma) = 41$, one of the two factors must have odd dimension. The only odd-dimensional irreducible symmetric space of noncompact type and rank $2$ is $SL_3(\RR)/SO_3$. Using again Proposition \ref{hyperbolicfactors} we see that there exists a reflective submanifold $\Sigma'$ of $M$ with $\dim(\Sigma') \geq \dim(\Sigma)$, which is a contradiction.

Consequently $\Sigma$ cannot exist and it follows that $i(M) = 24 = i_r(M)$.
\end{proof}

\begin{thm} \label{E7E6U1}
For $M = E_7^{-25}/E_6U_1$ we have $i(M) = i_r(M) = 22$.
\end{thm}

\begin{proof}
We already know from \cite{BO16} that $i_r(M) = 22$. From Corollary \ref{intermediate} we know that $i(M) \geq 13$. Let $\Sigma$ be a maximal totally geodesic submanifold of $M$ and assume that $\codim(\Sigma) \in \{13,\ldots,21\}$, that is $\dim(\Sigma) \in \{33,\ldots,41\}$. From the results in \cite{BO16} we can assume that $\Sigma$ is semisimple. Using Corollary \ref{intermediate} we obtain $(\rk(M) - \rk(\Sigma)) + \ell_\Sigma = \Lambda^M_\Sigma  \leq 16$. Since $\rk(M) = 3$, we have $\ell_\Sigma \leq 14$ if $\rk(\Sigma) = 1$, $\ell_\Sigma \leq 15$ if $\rk(\Sigma) = 2$, and $\ell_\Sigma \leq 16$ if $\rk(\Sigma) = 3$. 

The hyperbolic spaces $\FF H^k$ with $\ell_{\FF H^k} \leq 16$ are $\RR H^k$ for $k \in \{2,3,4,5,6,7\}$ (then $\ell_{\RR H^k} \in \{0,1,3,6,10,15\}$), $\CC H^k$ for $k \in \{2,3,4,5\}$ (then $\ell_{\CC H^k} \in \{1,4,9,16\}$), and $\HH H^k$ for $k \in \{2,3\}$ (then $\ell_{\HH H^k} \in \{6,13\}$). Since $\dim(\Sigma) \geq 33$ and $\ell_\Sigma \leq 16$, we easily see that $\Sigma$ cannot have rank $1$ or be a Riemannian product of rank $1$ symmetric spaces.

Assume that $\rk(\Sigma) = 2$. Since $\Sigma$ cannot be a Riemannian product of two rank $1$ symmetric spaces, $\Sigma$ is irreducible. The irreducible rank $2$ symmetric spaces $\Sigma$ with $\dim(\Sigma) \in \{33,\ldots,41\}$ are
$SO^o_{2,q}/SO_2SO_q$ ($q \in \{17,18,19,20\}$), $SU_{2,q}/S(U_2U_q)$ ($q \in \{9,10\}$) and $Sp_{2,5}/Sp_2Sp_5$. In all cases we have $\ell_\Sigma > 15$ (see Table \ref{symspacedata}) and therefore no such $\Sigma$ exists.  

Assume that $\rk(\Sigma) = 3$. We already know that $\Sigma$ cannot be the product of three rank $1$ symmetric spaces. Assume that $\Sigma = \Sigma_1 \times \Sigma_2$, where $\Sigma_1$ has rank $1$ and $\Sigma_2$ has rank $2$. Using Proposition \ref{hyperbolicfactors} we see that $\Sigma_1$ must be $\HH H^2$ or $\HH H^3$. Thus $\dim(\Sigma_1) \in \{8,12\}$ and hence $21 \leq \dim(\Sigma_2) \leq 33$ and $\ell_{\Sigma_2} \leq 10$. Using Table \ref{symspacedata} we see that the only possibility is $\Sigma_2 = Sp_{2,3}/Sp_2Sp_3$, which satisfies $\dim(\Sigma_2) = 24$ and $\ell_{\Sigma_2} = 9$. However, since $\dim(\Sigma) \geq 33$, we must have $\Sigma_1 = \HH H^3$, which gives $\ell_\Sigma = \ell_{\Sigma_1} + \ell_{\Sigma_2} = 13 + 9 = 22$, which is a contradiction. We conclude that $\Sigma$ cannot be the Riemannian product of a rank $1$ symmetric space and a rank $2$ symmetric space. Finally, assume that $\Sigma$ is irreducible. Using Table \ref{symspacedata} we obtain that the irreducible rank $3$ symmetric spaces $\Sigma$ with $\dim(\Sigma) \in \{33,\ldots,41\}$ and $\ell_\Sigma \leq 16$ are $\Sigma = SU_{3,6}/S(U_3U_6)$ (then $\dim(\Sigma) = 36$ and $\ell_\Sigma = 11$) and $\Sigma = Sp_{3,3}/Sp_3Sp_3$ (then $\dim(\Sigma) = 36$ and $\ell_\Sigma = 9$). We have $\Omega_M = 26$. For $\Sigma = SU_{3,6}/S(U_3U_6)$ we have $\Lambda_M^\Sigma = \ell_\Sigma = 11$ and hence $\codim(\Sigma) = 18 < 18.5 = \frac{1}{2}(\Omega_M + \Lambda_M^\Sigma)$, which contradicts Proposition \ref{firstestimate}. Thus $\Sigma = SU_{3,6}/S(U_3U_6)$ is not possible. The restricted root system of $M = E_7^{-25}/E_6U_1$ is of type ($C_3$), the six short roots have multiplicity $8$ and the three long roots have multiplicity $1$. The restricted root system of $\Sigma = Sp_{3,3}/Sp_3Sp_3$ is also of type ($C_3$), but the six short roots have multiplicity $4$ and the three long roots have multiplicity $3$. Due to the multiplicities of the long roots, it is not possible to realize the second root system as a subsystem of the first one, which implies that $\Sigma = Sp_{3,3}/Sp_3Sp_3$ is not possible either.
\end{proof}

\begin{thm} \label{E8E7Sp1}
For $M = E_8^{-24}/E_7Sp_1$ we have $i(M) = i_r(M) = 48$.
\end{thm}

\begin{proof}
We already know from \cite{BO16} that $i_r(M) = 48$. From Corollary \ref{intermediate} we know that $i(M) \geq 42$. Let $\Sigma$ be a maximal totally geodesic submanifold of $M$ and assume that $\codim(\Sigma) \in \{42,\ldots,47\}$, that is, $\dim(\Sigma) \in \{65,\ldots,70\}$. From the results in \cite{BO16} we can assume that $\Sigma$ is semisimple. Using Corollary \ref{intermediate} we obtain $(\rk(M) - \rk(\Sigma)) + \ell_\Sigma = \Lambda^M_\Sigma  \leq 10$. Since $\rk(M) = 4$, we have $\ell_\Sigma \leq 7$ if $\rk(\Sigma) = 1$, $\ell_\Sigma \leq 8$ if $\rk(\Sigma) = 2$, $\ell_\Sigma \leq 9$ if $\rk(\Sigma) = 3$ and $\ell_\Sigma \leq 10$ if $\rk(\Sigma) = 4$. 

The hyperbolic spaces $\FF H^k$ with $\ell_{\FF H^k} \leq 10$ are $\RR H^k$ for $k \in \{2,3,4,5,6\}$ (then $\ell_{\RR H^k} \in \{0,1,3,6,10\}$), $\CC H^k$ for $k \in \{2,3,4\}$ (then $\ell_{\CC H^k} \in \{1,4,9\}$), and $\HH H^2$ (then $\ell_{\HH H^2} = 6$). Since $\dim(\Sigma) \geq 65$ and $\ell_\Sigma \leq 10$, we easily see that $\Sigma$ cannot have rank $1$ or be a Riemannian product of rank $1$ symmetric spaces.

Assume that $\rk(\Sigma) = 2$. Since $\Sigma$ cannot be a Riemannian product of two rank $1$ symmetric spaces, $\Sigma$ is irreducible. The irreducible rank $2$ symmetric spaces $\Sigma$ with $\dim(\Sigma) \in \{65,\ldots,70\}$ are
$SO^o_{2,q}/SO_2SO_q$ ($q \in \{33,34,35\}$) and $SU_{2,17}/S(U_2U_{17})$. In all cases we have $\ell_\Sigma > 8$ (see Table \ref{symspacedata}) and therefore no such $\Sigma$ exists. 

Assume that $\rk(\Sigma) = 3$. We already know that $\Sigma$ cannot be the product of three rank $1$ symmetric spaces. Assume that $\Sigma = \Sigma_1 \times \Sigma_2$, where $\Sigma_1$ has rank $1$ and $\Sigma_2$ has rank $2$. We must have $\ell_{\Sigma_1} + \ell_{\Sigma_2} = \ell_\Sigma \leq 9$. Using Proposition \ref{hyperbolicfactors} we obtain that $\Sigma_1 = \HH H^2$, which implies $\dim(\Sigma_2) \in \{57,\ldots,62\}$ and $\ell_{\Sigma_2} \leq 3$. From Table \ref{symspacedata} we see that no such $\Sigma$ exists. Finally, for the case that $\Sigma$ is irreducible, we see from Table \ref{symspacedata}  that there exists no such $\Sigma$ with $\dim(\Sigma) \in \{65,\ldots,70\}$ and $\ell_\Sigma \leq 9$.

Assume that $\rk(\Sigma) = 4$. We already know that $\Sigma$ cannot be the product of four rank $1$ symmetric spaces. Assume that $\Sigma = \Sigma_1 \times \Sigma_2 \times \Sigma_3$, where $\Sigma_1,\Sigma_2$ have rank $1$ and $\Sigma_3$ has rank $2$. We must have $\ell_{\Sigma_1} + \ell_{\Sigma_2} + \ell_{\Sigma_3} = \ell_\Sigma \leq 10$. It follows that at least one of the two rank $1$ symmetric spaces is a real or complex hyperbolic space. Proposition \ref{hyperbolicfactors} then implies that this case cannot occur.  Assume that $\Sigma = \Sigma_1 \times \Sigma_2$, where $\Sigma_1$ has rank $1$ and $\Sigma_2$ has rank $3$. Using Proposition \ref{hyperbolicfactors} we obtain that $\Sigma_1 = \HH H^2$, which implies $\dim(\Sigma_2) \in \{57,\ldots,62\}$ and $\ell_{\Sigma_2} \leq 4$. From Table \ref{symspacedata} we see that no such $\Sigma$ exists. Assume that $\Sigma$ is irreducible. Then $\dim(\Sigma) \in \{65,\ldots,70\}$, $\rk(\Sigma) = 4$ and $\ell_\Sigma \leq 10$, and from Table \ref{symspacedata} we see that no such $\Sigma$ exists. Finally, assume that $\Sigma = \Sigma_1 \times \Sigma_2$, where $\Sigma_1$ and $\Sigma_2$ are irreducible and of rank $2$. Then one of the two spaces, say $\Sigma_1$, must satisfy $\dim(\Sigma_1) \geq 33$, $\rk(\Sigma_1) = 2$ and $\ell_{\Sigma_1} \leq 10$. From Table \ref{symspacedata} we see that no such $\Sigma$ exists.
\end{proof}

This finishes the proof of Theorem \ref{mainexc}.


\begin{thebibliography}{[99]}

\bibitem{BCO16} 
J.\ Berndt, S.\ Console, C.\ Olmos: 
\emph{Submanifolds and holonomy. Second edition}. 
Monographs and Research Notes in Mathematics. CRC Press, Boca Raton, FL, 2016.

\bibitem{BO16}
J.~Berndt, C.~Olmos:  
Maximal totally geodesic submanifolds and index of symmetric spaces. 
\emph{J.\ Differential Geom.} \textbf{104} (2016), no.\ 2, 187--217.

\bibitem{BO17}
J.~Berndt, C.~Olmos:  
The index of compact simple Lie groups.
\emph{Bull.\ London Math.\ Soc.} \textbf{49} (2017), 903--907.

\bibitem{BO18}
J.~Berndt, C.~Olmos:  
On the index of symmetric spaces.
\emph{J.\ Reine Angew.\ Math.} \textbf{737} (2018), 33--48.

\bibitem{CN78}
B.Y.~Chen, T.~Nagano:  
Totally geodesic submanifolds of symmetric spaces. II. 
\emph{Duke Math.\ J.} \textbf{45} (1978), 405--425.

\bibitem{E17} 
J.~Eschenburg:
Extrinsic symmetric spaces.
Preprint, 2017, \url{http://myweb.rz.uni-augsburg.de/~eschenbu/extsym.pdf}.

\bibitem{H01}
S.~Helgason: 
\emph{Differential geometry, Lie groups, and symmetric spaces}. 
Corrected reprint of the 1978 original. Graduate Studies in Mathematics, 34. American
Mathematical Society, Providence, RI, 2001.

\bibitem{IT00}
O.~Ikawa, H.~Tasaki:
Totally geodesic submanifolds of maximal rank in symmetric spaces.
\emph{Japan.\ J.\ Math.} \textbf{26} (2000), 1--29.

\bibitem{Iw66}
N.~Iwahori: 
On discrete reflection groups on symmetric Riemannian manifolds, in:
\emph{Proc.\ US-Japan Seminar in Differential Geometry (Kyoto 1965)}, Nippon Hyoronsha, Tokyo (1966), 57--62.

\bibitem{K08}
S.~Klein:
Totally geodesic submanifolds of the complex quadric.
\emph{Differential Geom.\ Appl.} \textbf{26} (2008), no.\ 1, 79--96.

\bibitem{K09}
S.~Klein: 
Reconstructing the geometric structure of a Riemannian symmetric space from its Satake diagram.
\emph{Geom.\ Dedicata} \textbf{138} (2009), 25--50.

\bibitem{K10a}
S.~Klein: 
Totally geodesic submanifolds of the complex and the quaternionic $2$-Grass\-mannians.
\emph{Trans.\ Amer.\ Math.\ Soc.} \textbf{361} (2010), no.\ 9, 4927--4967.

\bibitem{K10b}
S.~Klein: Totally geodesic submanifolds of the exceptional Riemannian symmetric spaces of rank $2$.
\emph{Osaka J.\ Math.} {\bf 47} (2010), no.\ 4, 1077--1157.

\bibitem{KN64}
S.~Kobayashi, T.~Nagano:
On filtered Lie algebras and geometric structures. I. 
\emph{J.\ Math.\ Mech.} \textbf{13} (1964), 875--907. 

\bibitem{L75} D.S.P.~Leung:
On the classification of reflective submanifolds of Riemannian symmetric spaces.
\emph{Indiana Univ.\ Math.\ J.} \textbf{24} (1974/75), 327--339. 
Errata: \emph{Indiana Univ.\ Math.\ J.} \textbf{24} (1975), no.\ 12, 1199.

\bibitem{L79} D.S.P.~Leung:
Reflective submanifolds. III. Congruency of isometric reflective submanifolds and corrigenda to the classification of reflective submanifolds. 
\emph{J.\ Differential Geom.} \textbf{14} (1979), no.\ 2, 167--177.

\bibitem{NT95} 
T.~Nagano, M.S.~Tanaka: 
The involutions of compact symmetric spaces. III.
\emph{Tokyo J.\ Math.} \textbf{18} (1995), 193--212. 

\bibitem{O80}
\textcyr{A.L.~Onishchik: 
O vpolne geodezicheskih podmnogoobraziyah simmetricheskih prostranstv.
\emph{Geometricheskie metody v zadachah analiza i algebry}} \textbf{2} (1980), 64--85.

\bibitem{S62} 
J.~Simons:
On the transitivity of holonomy systems.  
\emph{Ann.\ of Math.\ (2)} \textbf{76} (1962), 213--234.

\bibitem{Wo63} J.A.~Wolf:
Elliptic spaces in Grassmann manifolds.
\emph{Illinois J.\ Math.} \textbf{7} (1963), 447--462. 

\end{thebibliography}
\end {document}